\documentclass[11pt]{article}

\usepackage{cite}
\usepackage{amssymb,amsfonts,amsmath,amsthm}
\usepackage{geometry}
\usepackage{color,enumerate,graphicx,hyperref}
\usepackage{stix}
\usepackage[normalem]{ulem}
\usepackage{tikz}
\usetikzlibrary{hobby}

\newcommand{\eps}{\varepsilon}
\newcommand{\R}{\mathbb R}

\newcommand{\Sphere}{\mathbb S}
\newcommand{\Z}{\mathbb Z}

\newcommand{\ee}{\lambda}

\newcommand{\x}{x}
\newcommand{\m}{m}

\newcommand{\mpe}{m'}

\newcommand{\intd}{\, \mathrm{d}}

\newtheorem{theorem}{Theorem}
\newtheorem{proposition}[theorem]{Proposition}
\newtheorem{definition}[theorem]{Definition}
\newtheorem{lemma}[theorem]{Lemma}
\newtheorem{corollary}[theorem]{Corollary}
\newtheorem{remark}[theorem]{Remark}

\numberwithin{equation}{section}
\numberwithin{theorem}{section}

\usepackage{todonotes}

\title{Uniqueness and radial symmetry of magnetic skyrmions}

\author{Giovanni Di Fratta\footnote{Dipartimento di Matematica e
    Applicazioni, “R. Caccioppoli”, Universit\`a degli Studi di Napoli
    “Federico II”, Via Cintia, 80126 Napoli, Italy} \and Cyrill
  B. Muratov\footnote{Dipartimento di Matematica, Universit\'a di
    Pisa, Largo Bruno Pontecorvo 5, 56127 Pisa, Italy} \and Theresa
  M. Simon\footnote{Institut f\"{u}r Analysis und Numerik,
    Universit\"{a}t M{\"u}nster, 48149 M{\"u}nster, Germany} \and
  Valeriy V. Slastikov\footnote{School of Mathematics, University of
    Bristol, Bristol BS8 1TW, United Kingdom}}

\date{\today}

\begin{document}
\maketitle

\begin{abstract}
  This paper establishes radial symmetry and uniqueness of magnetic
  skyrmions as global minimizers of the simplest local micromagnetic
  energy. The strategy is to show that any radial critical point with
  strictly decreasing out-of-plane component is the unique minimizer
  in the class of general degree-one maps with the prescribed limit at
  infinity. To this end, we identify a class of null Lagrangians that
  is sufficiently flexible to contain a calibration, giving a uniform
  lower bound for the energy that is uniquely saturated on the
  critical point. Our arguments provide a blueprint for proofs of
  symmetry and uniqueness of minimizers in vectorial variational
  problems under topological constraints.
\end{abstract}



\tableofcontents

\section{Introduction}
\label{sec:introduction}

Coherent structures are a central object of interest in a great many
problems described by nonlinear partial differential equations (PDEs)
and are key to understanding the behavior of the physical systems
modeled by these equations
\cite{sulem,tinkham,braun12,jaffe,majda,cross93,ko:book}. They appear
as particular solutions of these PDEs that exhibit varying degrees of
spatial localization. Examples of coherent structures arising in the
models of three-dimensional physical media include transition layer
solutions approximating two-dimensional surfaces of discontinuity
between two distinct phases in the phase field models of phase
transitions \cite{modica87}, vortex solutions of Ginzburg--Landau
theory representing vortex filaments in type-II superconductors
\cite{bethuel,baldo12}, or ground state (bright soliton) solutions of
the nonlinear Schr\"odinger (NLS) equation, representing point-like
concentrations of energy \cite{tao09}.

The existence of such localized solutions typically requires a
delicate balance between different terms in the equation, such as
dispersion and the nonlinearity for the NLS solitons. As a
consequence, these solutions are often fragile: a small change in the
parameters of the model or in the initial data may disturb this
balance and lead to the destruction of the coherent structure, which
may then disperse, collapse, or relax to a spatially uniform
state. For instance, the NLS ground states are unstable in the
mass-supercritical regime, where arbitrarily small perturbations may
lead to finite-time blow-up or dispersion \cite{nakanishi11}, while a
closed vortex filament in a type-II superconductor may shrink and
disappear under the action of its line tension \cite{bethuel:ann06}. A
distinguished class of coherent structures that do not share this
fragility is formed by {\em topological solitons}
\cite{jaffe,manton}. These are localized solutions of nonlinear PDEs
whose unknowns take values in a manifold with non-trivial topology and
are often global minimizers of an energy functional over the
respective topological sectors. The latter may be characterized by a
non-zero integer topological invariant, such as the winding number or
the degree of the map describing a configuration. Since this invariant
cannot be altered by continuous deformations of the configuration, a
topological soliton cannot be continuously deformed into the
topologically trivial state, and its existence and stability are to a
large extent insensitive to the details of the model. Moreover, the
energy in each non-trivial topological sector is often bounded from
below by a multiple of the absolute value of the topological charge,
so that a topological soliton cannot decay without either overcoming
an energy barrier or passing through a singular configuration.

One of the first known examples of topological solitons is given by
skyrmions, introduced by Skyrme as minimizers of an energy functional
defined on maps from $\mathbb{R}^3$ to $\mathbb{S}^3$ that converge to
a constant at infinity and, hence, carry an integer-valued degree
\cite{skyrme61,skyrme62,esteban86,esteban90,esteban04,lin04a,manton}. Their
two-dimensional analogs are maps $\m : \mathbb{R}^2 \to \mathbb{S}^2$
approaching a fixed point of $\mathbb{S}^2$ at infinity, whose
topological charge, for smooth maps or more generally for maps with
finite Dirichlet energy, is the Brouwer degree, which admits the
integral representation \cite{brezis95}
\begin{align}
  \label{eq:N}
  \mathcal{N}(\m) = \frac{1}{4\pi} \int_{\mathbb{R}^2}
  \m \cdot (\partial_1 \m \times \partial_2
  \m) \, dx \in \mathbb{Z}.
\end{align}
The Dirichlet energy of such maps satisfies the lower bound
$\int_{\mathbb{R}^2} |\nabla \m|^2 \, d x \geq 8 \pi
|\mathcal{N}(\m)|$, with equality if and only if $\m$ is conformal or
anti-conformal, i.e., after the stereographic projection, a rational
or anti-rational function of degree $|\mathcal{N}(\m)|$
\cite{belavin75,eells78,brezis85}. Since the Dirichlet energy is
conformally invariant in two dimensions, these minimizers form a
non-compact family with an arbitrary scale. By Derrick's scaling
argument \cite{derrick64}, adding to the energy a lower order term
that is minimized by the uniform state, such as an anisotropy term
(see below), makes it energetically favorable for these maps to
collapse, so that non-trivial minimizers fail to exist. Stable
two-dimensional skyrmions thus require an additional term in the
energy that selects their scale, such as the higher order term of the
baby Skyrme model \cite{lin04,li11}.

Alternatively, topologically non-trivial maps
$\m : \R^2 \to \mathbb S^2$ may be stabilized by a bilinear term in
$\m$ and $\nabla \m$ that naturally arises in micromagnetics
\cite{bogdanov89,bogdanov89a}. This so-called Dzyaloshinskii--Moriya
interaction (DMI) \cite{dzyaloshinskii58,moriya60} breaks the
reflection symmetry of the energy, resulting in the emergence of
topologically non-trivial point-like configurations in the plane
\cite{melcher14,li18, greco19,bms:prb20,gustafson21,bms:arma21} that
are referred to as {\em magnetic skyrmions}. Because of their
topological protection, magnetic skyrmions have received a great deal
of attention, both theoretically and experimentally, as promising
candidates for information technology applications
\cite{nagaosa13,fert13,fert17,zhang20}.  In the micromagnetic context,
$\m$ is the normalized magnetization vector, and in a skyrmion
configuration it minimizes the associated micromagnetic energy
functional among all configurations with topological degree
$\mathcal N(\m) = 1$, with the convention $\m(\infty) = -e_3$ to fix
the sign of the degree. Specifically, the simplest model appropriate
for ultrathin ferromagnetic films with perpendicular
magnetocrystalline anisotropy and the DMI of interfacial origin is
described by the following energy \cite{rohart13,ms:prsa17,bms:prb20}:
\begin{align}
  \label{eq:E}
  \widetilde E_{\kappa,Q}(\m) = \int_{\mathbb{R}^2} \left( |\nabla \m|^2 +
  (Q - 1) |\m'|^2 - 2 \kappa \m' \cdot
  \nabla m_3 \right) d x,
\end{align}
where $\m = (\m', m_3)$, with $\m' = (m_1, m_2) \in \mathbb{R}^2$ and
$m_3 \in \mathbb{R}$, $Q > 1$ is the ratio of the magnetocrystalline
and shape anisotropy constants, and $\kappa > 0$ (without loss of
generality) measures the DMI strength.

Existence of minimizers among all $\m : \R^2 \to \mathbb S^2$ such
that $\m + e_3 \in H^1(\R^2; \R^3)$ and $\mathcal N(\m) = 1$ was
established in \cite{bms:prb20,bms:arma21} for all $Q > 1$ and
$0 < \kappa < \sqrt{(Q-1)/2}$, and in the maximal range
$0 < \kappa < {4 \over \pi} \sqrt{Q - 1}$ among axisymmetric profiles
in \cite{dmss:rad}. Furthermore, it was shown in \cite{bms:arma21}
that in the {\em conformal limit} $\kappa \to 0$ these minimizers
concentrate on points and become asymptotically radial, confirming the
physical intuition \cite{rohart13}. In a closely related setting where
$\m \in H^1(\Omega; \mathbb S^2)$ with $\Omega \subset \R^2$ a bounded
domain with sufficiently smooth boundary (e.g., a disk), and with
Dirichlet boundary condition $\m = -e_3$ on $\partial \Omega$,
minimizers of degree one were shown to concentrate, as $\kappa \to 0$
and $Q \to 1$ jointly with a suitable rate, on a single point
$x \in \Omega$, around which the profile becomes asymptotically radial
\cite{mmss:cmp23}. Moreover, as $Q \to \infty$ for $\kappa > 0$ and
all positive degrees $d$ below a certain threshold, the energy density
was shown \cite{mss:arma25} to concentrate on a weighted sum of Dirac
deltas $8 \pi \sum_{i=1}^k d_i \delta_{x_i}$, for a finite set of
points $x_i \in \overline\Omega$ and multiplicities
$d_i \in \mathbb N$ obeying $\sum_{i=1}^k d_i = d$. For
$x \in \Omega$, the term $8 \pi \delta_{x}$ is the limiting energy
density of a single Belavin--Polyakov profile (degree one harmonic map)
concentrating at $x$.

Whether magnetic skyrmion minimizers in the whole plane are indeed
radially symmetric or whether their profile is unique has been, to the
best of our knowledge, a difficult open problem. It is an example of
an important question that is relevant across problems exhibiting
coherent structures, and precious few results of this type currently
exist. Here one should mention celebrated results on radial symmetry
of ground states for NLS type problems \cite{gidas79,gidas81}
(uniqueness among radial profiles is still another issue; for fairly
general results, see \cite{serrin00,lewin20}), which rely on the
method of moving planes. In the vectorial setting, where this
technique is generally not available beyond cooperative systems
\cite{troy81}, results of this kind are much scarcer.

A notable exception is the Ginzburg--Landau equation in two
dimensions, whose equivariant vortex solutions of degree
$d \in \mathbb{Z}$ at infinity are uniquely determined by the
associated ordinary differential equation in polar coordinates
\cite{herve94,chen94}. Building on the quantization result of Brezis,
Merle and Rivi\`ere \cite{brezis94}, Mironescu proved that every local
minimizer (in the sense of De Giorgi) with degree $d = \pm 1$ of the
Ginzburg--Landau energy coincides, up to a translation, a constant
phase and complex conjugation, with the radial vortex
\cite{mironescu96} (see also \cite{delpino04}), while Sandier showed
that non-constant local minimizers necessarily have degree $d = \pm 1$
\cite{sandier98}, yielding a complete classification. On the disk
$B_R \subset \mathbb{R}^2$ of radius $R > 0$ with the degree-one
boundary condition $u(x) = x/|x|$ for $x \in \partial B_R$, Matano
asked whether the equivariant solution is the unique minimizer of the
Ginzburg--Landau energy for every $R > 0$
\cite[Problem~10]{bethuel}. This was known for $R$ sufficiently small,
by convexity \cite[Section~2]{brezis23}, and for $R$ sufficiently
large the answer was obtained in \cite[Theorem~10.2]{pacard00}, but
the full answer, in the affirmative, was provided only very recently
in the preprint \cite{chen26}, whose authors acknowledge the use of AI
tools. Other notable results in higher dimensions include
\cite{millot10,ignat18a}, as well as a very recent preprint
\cite{ignat26}. Related symmetry and classification results for
minimizers of micromagnetically motivated Dirichlet-type energies have
been obtained on spherical and cylindrical surfaces
\cite{difratta19,digratta23}, on surfaces of revolution under
additional average conditions \cite{difratta24}, and on bounded
domains with boundary penalization \cite{difratta22}.

The goal of the present paper is to investigate the symmetry and
uniqueness of magnetic skyrmions in the plane, a problem that is
closely analogous to that for two-dimensional Ginzburg--Landau
vortices discussed above. The radial symmetry is understood as the
special (corotational) form of the energy-minimizing profile (see the
next section for further details):
\begin{align}
  \label{eq:equivar}
  \m_*(x) = -\left( {x \over |x|} \sin\theta_* (|x|) , \cos
    \theta_*(|x|) \right) \qquad x \in \R^2,
\end{align}
for some angle function $\theta_* : [0, \infty] \to [0, \pi]$ that is
monotonically decreasing from $\theta_*(0) = \pi$ to
$\theta_*(\infty) = 0$, with the convention $\m_*(0) = e_3$.  We give
a positive answer to the uniqueness and radial symmetry question for
the problem in the whole plane in the following theorem. An analogous
positive answer for the problem on a disk for a range of parameters
will be presented elsewhere \cite{dmss:rad2}.

\begin{theorem}
  \label{thm:main}
  Let $Q > 1$ and let $0 < \kappa < {4 \over \pi} \sqrt{Q - 1}$. Then
  the energy $\widetilde E_{\kappa,Q}(\m)$ in \eqref{eq:E} admits a
  minimizer $\m_*$ among all $m : \R^2 \to \mathbb S^2$ with
  $\m + e_3 \in H^1(\R^2; \R^3)$ and $\mathcal N(\m) =
  1$. Furthermore, the minimizer $\m_*$ is unique, up to translations,
  and has the form \eqref{eq:equivar}, where
  $\theta_* \in C^\infty([0,\infty))$ satisfies $\theta_*(0) = \pi$,
  $\theta_*(\infty) = 0$ and $\theta_*' < 0$ on $[0,\infty)$. In
  particular, $m_* \in C^\infty(\R^2; \mathbb S^2)$.
\end{theorem}

\noindent We note that the range of $\kappa$ is sharp, as the energy
is unbounded below for $\kappa > {4 \over \pi} \sqrt{Q - 1}$ (see
\cite{ms:prsa17,dmss:rad}) and the infimum of the energy is zero and
is not attained for $\kappa = {4 \over \pi} \sqrt{Q - 1}$ (see
Proposition \ref{prop:crit}).

The proof of Theorem \ref{thm:main} is obtained by constructing a
calibration of the energy in \eqref{eq:E} in terms of the minimizer of
this energy among all radially symmetric degree one profiles, whose
existence for not too big values of $\kappa / \sqrt{Q - 1}$ was
established in \cite{gustafson21}, and in the full range of $\kappa$
was established in \cite{dmss:rad}. The analysis of the associated
Euler--Lagrange equation reveals monotonicity of the energy minimizing
angle profile (see also \cite{gustafson21}), which allows to build a
calibration via a pointwise lower bound on the energy density, using
expressions that involve null Lagrangians associated with the
topological invariants of the considered problem. The fact that the
null Lagrangian is indeed a lower bound then reduces to certain
inequalities on the radially symmetric profile $\theta_*$, which are
amenable to standard ODE arguments.

\begin{remark}
  \label{rem:equivar}
  We note that in contrast to the Ginzburg-Landau equation, minimizers
  of \eqref{eq:E} with non-zero degree $d = \mathcal N(m) \not= 1$
  cannot be captured by passing to the respective equivariant
  profiles. An equivariant profile can be written in the (inverted)
  spherical coordinates as
  \begin{align}
    \label{eq:equivarsph}
    m = -(\sin \theta \cos \phi, \sin \theta \sin \phi, \cos \theta),    
  \end{align}
  where $\theta = \theta (r)$ and $\phi = \psi d + \phi_0$ are the
  polar and azimuthal angle, respectively, as functions of the polar
  coordinates $r = |x + i y|$ and $\psi = \arg(x + i y)$ in the plane,
  with $\theta(0) = \pi$, $\theta(\infty) = 0$ and
  $\phi_0 \in \mathbb R$. However, it can be readily verified that the
  Euler-Lagrange equation for the energy in \eqref{eq:E} is not
  satisfied by this ansatz unless $d = 1$ and
  $\phi_0 \in \pi \mathbb Z$. Therefore, minimizers with higher
  degrees, if they exist at all in the whole plane, cannot be
  described by radial (equivariant) profiles (see also \cite{li18}).
\end{remark}

\paragraph{Acknowledgements} \textsc{G.D.F.} is a member of
GNAMPA--INdAM. He gratefully acknowledges partial financial support
from the GNAMPA Project CUP\_E53C25002010001, and from the University
of Naples Federico II through the FRA Project-B ``VarMoCry'' on
\emph{Variational Analysis and Modeling of Liquid Crystals}.
C.B.M. is a member of INdAM/GNAMPA and acknowledges the MUR Excellence
Department Project awarded to the Department of Mathematics,
University of Pisa, CUP I57G22000700001. T.M.S. is funded by the
Deutsche Forschungsgemeinschaft (DFG, German Research Foundation)
under Germany's Excellence Strategy EXC 2044/2 –390685587, Mathematics
Münster: Dynamics–Geometry–Structure. C.B.M. is grateful to F. Rybakov
for valuable discussions and to A. Bernand-Mantel for providing the
simulation results in Fig. \ref{fig:mins}.  The authors acknowledge
the use of AI tools (Claude Opus 5 and ChatGPT-6 Astra) to assist in
verifying proofs and simplifying technicalities that emerged when
testing specific approaches. All mathematical ideas and proof
strategies were developed solely by the authors, who bear full
responsibility for the content of this paper.

\section{Setup and strategy}\label{sec:main_results}

Since in this paper we consider configurations defined on the entire
plane $\R^2$, we can always set $Q=2$ by a rescaling, in which case
the DMI constant becomes $\bar \kappa := \frac{\kappa}{\sqrt{Q-1}}$
(see, e.g., \cite{dsbm:jns26}).  Hence for $\bar\kappa>0$ and
$\m : \R^2 \to \mathbb S^2$ we define the energy
\begin{align}\label{eq:Energy}
  E_{\bar \kappa}(m) := 
  \int_{\mathbb{R}^2} \left( |\nabla \m|^2 + 
  |\m'|^2 - 2 \bar \kappa \m' \cdot
  \nabla m_3 \right) \intd x,
\end{align}
and the associated energy density
\begin{equation}\label{eq:energy-density}
 e_{\bar \kappa}(\m):=|\nabla \m|^2+|\mpe|^2-2\bar \kappa
 \mpe\cdot\nabla m_3. 
\end{equation}
Here and everywhere below we always mean
$v \cdot \nabla = \sum_{i=1}^2 v_i \partial_i$ for a vector field
$v : \R^2 \to \R^2$. This energy is well defined over the admissible
class
\begin{align}\label{eq:admissible}
 \mathcal A:=\left\{\m\in \mathring{H}^1(\R^2;\mathbb{S}^2): m + e_3
  \in L^2(\R^2;\R^3)\right\},
\end{align}
and the corresponding degree-one class is
\begin{align}\label{eq:degree-class}
 \mathcal A_1& :=\{\m\in\mathcal A:\mathcal N(\m)=1\}.
\end{align}
We define {\em radial} magnetizations as
\begin{align}
	\label{eq:equivar_general}
	\m_\theta(x) = -\left( \sin\theta (|x|) \frac{x}{|x|} , \cos
	\theta(|x|) \right),
\end{align}
where $x\in \R^2$ and $\theta: [0, \infty) \to \R$. Strictly speaking,
these configurations belong to the class of corotational
configurations (a specific class of $1$-equivariant, or axially
symmetric, configurations), and in the physics literature they are
referred to as {\em N\'eel skyrmions}. Here we adopt the term ``radial
configurations'' just for simplicity (see also Remark
\ref{rem:equivar}).

The admissible class $\mathcal A$ for $\m$ of the form
\eqref{eq:equivar_general} is then
\begin{align}
  \label{eq:Ar} \mathcal{A}_\mathrm{r}:= \left\{ \theta
  \in H^1_\mathrm{loc}(0,\infty): r |\theta'|^2+ r (1 - \cos
  \theta) + \frac{\sin^2 
  \theta}{r}\in L^1(0,\infty)  
  \right\},
\end{align}
where we identify $\theta \in H^1_\mathrm{loc}(0,\infty)$ with its
continuous representative.  We define a subclass
\begin{align}
  \mathcal{A}_{\mathrm{r},1} := \left\{ \theta \in
  \mathcal{A}_\mathrm{r}: \lim_{r\to 0} \theta(r) = \pi, \ \lim_{r
  \to \infty} \theta(r) = 0 \right\}
\end{align}
that contains the corresponding profiles of degree one. This
definition fixes an arbitrary additive constant from $2 \pi \mathbb Z$
in the definition of $\theta$, and we note that
$\lim_{r \to 0} \theta(r)$ exists for any
$\theta \in \mathcal{A}_\mathrm{r}$ by Lemma \ref{lemma:degree}.  Thus
$\theta \in \mathcal{A}_{\mathrm{r},1}$ is continuous on
$[0, \infty]$, with the convention $\m_\theta(0) = e_3$ and
$\m_\theta(\infty) = -e_3$.  For $\theta \in \mathcal{A}_{\mathrm{r}}$
we also define the energy
$I_{\bar \kappa}(\theta) := E_{\bar \kappa}(m_\theta)$, with
$m_\theta$ from \eqref{eq:equivar_general}, or explicitly
\cite{bogdanov89,rohart13}
\begin{align}
  I_{\bar \kappa}(\theta) := 2\pi \int_0^\infty \left[ r |\theta'|^2 +
  \left(r +\frac{1}{r} 
  \right)\sin^2 \theta + 2 \bar \kappa r \theta' \sin^2\theta\right]  
  \intd r.
\end{align}

Our main result takes the surprisingly strong form: 
\begin{quote} \emph{
  If a strictly decreasing angle function
  $\theta_* \in \mathcal{A}_{\mathrm{r},1}$ is a \emph{critical point}
  of $I_{\bar \kappa}$, then $m_{\theta_*}$ is the unique minimizer of
  $E_{\bar \kappa}$ over $\mathcal{A}_1$, up to translations.}
\end{quote}

Note that, provided we have strict monotonicity of the profile
$\theta_*$, we do not need any assumptions on the parameters beyond
the purely qualitative requirement $\bar \kappa \in (0, \infty)$,
corresponding to $\kappa >0$ and $Q>1$ (in fact, up to a change
$m' \to -m'$ it is sufficient to assume $\kappa \not= 0$).  This
closely mimics the standard situation for strictly convex functionals,
where such a statement is true without any assumptions on the
structure of the critical point.

Standard regularity arguments show that every critical point of
$I_{\bar \kappa}$ from $\mathcal{A}_{\mathrm{r},1}$ is smooth
and solves
\cite{dmss:rad}
\begin{align}\label{eq:ELG}
  {d^2 \theta_* \over dr^2} + {1 \over r} {d
  \theta_* \over dr}  -  \left(1 + \frac{1}{r^2} \right)
  \sin\theta_* \cos \theta_* + {\bar \kappa \sin^2
  \theta_* \over r} = 0,
\end{align}
classically for every $r > 0$. This motivates the following
definition.
\begin{definition}\label{def:monotone_skyr}
  We say that \emph{$\bar \kappa >0$ supports a monotone skyrmion
    profile $\theta_*(r)$, if $\theta_* \in C^\infty(0, \infty)$} and
  the following conditions are satisfied:
  \begin{enumerate}[1)]    
  \item $\lim_{r \to 0} \theta_*(r) = \pi$ and
    $\lim_{r \to \infty} \theta_*(r) = 0$;
  \item $\theta'_*(r)<0$ for all $r > 0$;
  \item $\theta_*(r)$ satisfies \eqref{eq:ELG} for all $r > 0$.
  \end{enumerate}
\end{definition}

Monotone skyrmion profiles were first constructed numerically in
\cite{bogdanov89a} (see also \cite{rohart13,komineas23}).  Existence
of monotone skyrmion profiles for small $\bar\kappa > 0$ follows from
the results of Gustafson and Wang \cite{gustafson21}. In the full
range $\bar\kappa \in (0, 4/\pi)$, we recently established existence
of such profiles in the following theorem.
    
\begin{theorem}[Di Fratta, Muratov, Simon, Slastikov
  \cite{dmss:rad}]\label{thm:existence_mon_radial}
  Let $\bar \kappa >0$. Then $I_{\bar \kappa}$ admits a unique degree
  one minimizer $\theta_* \in \mathcal{A}_{\mathrm{r},1}$ if and only
  if $\bar \kappa < \frac{4}{\pi}$. The minimizer is of class
  $C^\infty([0, \infty))$, satisfies the Euler-Lagrange equation
  \eqref{eq:ELG}, as well as $\theta'_*(r) <0$ for all $r \geq
  0$. Furthermore, it is the unique monotone profile in the sense of
  Definition \ref{def:monotone_skyr}, and
  $m_{\theta_*} \in C^\infty(\R^2; \mathbb S^2)$.
\end{theorem}

\noindent Consequently, every $0 < \bar \kappa < \frac{4}{\pi}$
supports a unique monotone skyrmion profile $\theta_* (r)$. \medskip
    
The main result of this paper characterizes the unique global
minimizer as a monotone, radial skyrmion.
\begin{theorem}\label{thm:plane}
  Let $\bar \kappa>0$ support a monotone skyrmion profile
  $\theta_* (r)$. Then $m_{\theta_*}$ is the unique minimizer of
  $E_{\bar \kappa}$ over $\mathcal{A}_1$, up to translations.
\end{theorem}

\noindent As a corollary, we can immediately upgrade Theorem
\ref{thm:existence_mon_radial} to allow non-radial competitors, fully
solving the existence and symmetry question for degree one minimizers
of $E_{\bar \kappa}$ on the whole plane. This yields Theorem
\ref{thm:main}.

\begin{remark}\label{rem:topological-protection}
  The physical notion of topological protection of magnetic skyrmions
  stems from the fact that any global minimizer $m$ of
  $E_{\bar\kappa}$ in the degree one sector is automatically a
  \emph{local} minimizer among all $m \in \mathcal A$.  Indeed, by the
  continuity of the topological degree in the strong $H^1$-topology
  and its integer-valuedness, there exists $\delta>0$ such that every
  $m\in\mathcal A$ with $\|m-m_*\|_{H^1(\R^2)}<\delta$ still satisfies
  $\mathcal N(m)=1$.  Consequently,
  $E_{\bar\kappa}(m)\geqslant E_{\bar\kappa}(m_*)$. Thus, $m_*$ is
  energetically stable against small $H^1$-perturbations.
\end{remark}

\begin{remark}\label{rem:bulk-DMI}
  Theorem~\ref{thm:main} also applies to the two-dimensional models of
  magnetic materials with bulk DMI, in which the DMI energy density is
  proportional to $m \cdot (\nabla \times m)$ \cite{nagaosa13}. This
  energy is equivalent to the one considered in the present paper,
  modulo a 90-degree rotation of the magnetization vector (see the
  footnote in \cite{mss:arma25}).  As such a rotation fixes $-e_3$ and
  preserves the degree, Theorem~\ref{thm:main} yields the unique
  minimizer up to translations, of the corresponding energy. It has
  the same profile $\theta_*$ as in \eqref{eq:equivar}, with
  tangential in-plane magnetization (a Bloch skyrmion).
\end{remark}

The proof of Theorem \ref{thm:plane} is obtained by constructing a
calibration of the energy \eqref{eq:Energy} in terms of the monotone
skyrmion profile $\theta_* (r)$.  By this we mean a lower bound for
the energy that only involves null Lagrangians depending on the
(topological) data of our problem and the monotone profile
$\theta_*$. This lower bound has to be achievable on the radial
magnetization $m_{\theta_*} \in \mathcal A_1$.  We will describe it in
detail in Section \ref{sec:calibrations}.  This strategy reduces the
highly non-trivial problem of global minimality to studying the
properties of the monotone skyrmion profile $\theta_*$ using standard
ODE techniques.
\begin{figure}[t]
	\centering
	\includegraphics[width=14cm]{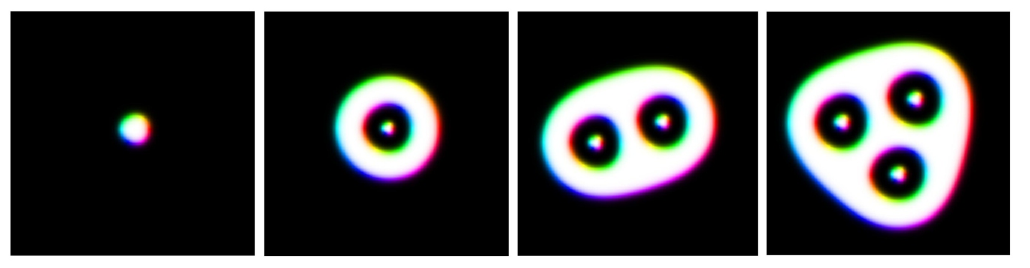}
	\caption{Several local minimizers of the energy in \eqref{eq:Energy}
		with $\bar\kappa = 1.2$ obtained by a numerical minimization,
		using {\sc Excalibur} software. The grayscale indicates the values
		of $m_3$, while the color code shows the direction of $m'$, with
		the first panel corresponding to the radial profile of the form
		\eqref{eq:equivar}. The computational domain shown is
		$57 \times 57$.} 
	\label{fig:mins}
\end{figure}

\begin{remark}
  While according to Theorem \ref{thm:plane} the unique monotone
  radial skyrmion profile determines the unique, up to translations,
  global minimizer of $E_{\bar \kappa}$ over $\mathcal A_1$, numerical
  evidence indicates that at least for some
  $0 < \bar \kappa < {4 / \pi}$ there are many non-radial {\em local}
  minimizers of $E_{\bar \kappa}$ in this class
  \cite{rybakov19,kuchkin21}. In Fig.~\ref{fig:mins}, we present a few
  such solutions, alongside with the canonical radial skyrmion
  solution obtained with the help of {\sc Excalibur} software
  \cite{excalibur}.
\end{remark}

\subsection{Calibrations and null Lagrangians}\label{sec:calibrations}

In this section we describe the basic idea behind calibrations and
provide useful definitions for the subsequent analysis. The
calibration method consists in constructing an auxiliary integrand
whose associated integral functional satisfies three properties:
\begin{enumerate}[\ \ 1)]
\item it provides a lower bound for the total energy of every
  admissible map;
\item it takes the same value on all maps in the admissible class;
\item it equals the energy of the expected minimizer.
\end{enumerate}
Together, these properties establish the global minimality of the
candidate. Crucially, the first two properties are required to hold
only after integration---once the boundary conditions and the
topological degree have been fixed---rather than pointwise at the
level of the integrand.
	
A natural starting point is therefore the theory of \emph{null
  Lagrangians}, whose integrals are left unchanged by compactly
supported variations. In the classical setting of Euclidean-valued
maps, first-order null Lagrangians depending only on the gradient are
precisely the affine combinations of its minors
\cite[Theorem~1]{olver88}.  In two space dimensions the gradient has
only two columns, so the only nonconstant minors are its entries and
its $2\times2$ determinants. The former are divergences, since
$\partial_j u_\alpha=\nabla\cdot(u_\alpha e_j)$, while the latter are
the Jacobian determinants of pairs of components.  This suggests
looking for calibrating integrands built from divergence terms and
Jacobians. Below we develop this construction for sphere-valued maps,
including Jacobians weighted by functions of the map, and verify that
their integrals are constant on the admissible class.

\begin{definition}
  We define the Jacobian determinant in the context of maps
  $\m\in \mathring H^1(\R^2;\Sphere^2)$ as
\begin{equation}\label{eq:degree}
  \omega_\m: =\m\cdot(\partial_1\m\times\partial_2\m).
\end{equation}
Then the degree \eqref{eq:N} is given by
\begin{align}
 \mathcal N(\m)=\frac1{4\pi}\int_{\R^2} \omega_\m\intd\x.
\end{align}
\end{definition}
Note that the Jacobian is integrable because almost everywhere on
$\R^2$ we have
\begin{align}\label{eq:jac-bound}
 |\omega_\m|\leqslant|\partial_1\m|\,|\partial_2\m|
 \leqslant\tfrac{1}{2}|\nabla \m|^2.
\end{align}
If we simply use (multiples of) $\omega_m$ as a calibration, this
gives the topological lower bound
\begin{align}\label{eq:top_bound}
	\int_{\R^2} |\nabla m|^2 \intd x \geq 8\pi|\mathcal{N}(m)|,
\end{align}
which Belavin and Polyakov used to characterize all minimizers of the
pure exchange energy model with a given degree, i.e., all harmonic
maps \cite{belavin75}.  However, this makes the construction we have
in mind far too rigid.

Instead, we observe that if $ \eta \in L^\infty(\mathbb{S}^2)$, then
the area formula (see for example \cite[Theorem 6.4]{mueller90} for a
version for Sobolev maps) implies
\begin{align}
  \int_{\R^2} \eta(m) \omega_m \intd x 
  = \int_{\mathbb{S}^2} \eta(z) q(z) \intd\mathcal{H}^2(z),
\end{align}
with $ q(z) := \sum_{m(x)=z} \operatorname{sgn}\, \omega_m(x) $.  It
turns out \cite[Lemma 4]{esteban92} that $q$ is constant almost
everywhere on $ \mathbb{S}^2 $ and equals the degree $\mathcal N(m)$
of the map $m$.  (There is the technical caveat that the paper
\cite{esteban92} is concerned with maps on $\R^3$, but the analogous
result also holds on $\R^2$; see also \cite{brezis99}.)  Therefore, we
have
\begin{align}
  \int_{\R^2} \eta(m) \omega_m \intd x = \mathcal{N}(m)
  \int_{\mathbb{S}^2} \eta(z) 
  \intd\mathcal{H}^2(z),
\end{align}
showing that $2 \eta(m) \omega_m$ is an admissible null Lagrangian,
where the factor of $2$ is included for convenience in later
computations.

The above observation significantly expands the class of possible
calibrations, but we still need to include the remaining option of
divergences.  On the whole space, constants cannot play any role
because we need our calibrations to have a reasonably defined
integral, although on bounded domains it is generally necessary to
include them \cite{dmss:rad2}.  In analogy with the observation above,
we should consider not only linear combinations of the $2 \times 2$
minors of $\nabla m$---which account for terms quadratic in the
gradient---but also divergence expressions of the form
$\nabla \cdot \Psi(m)$, with $\Psi : \mathbb{S}^2 \to \mathbb{R}^2$
sufficiently smooth, to capture terms linear in the gradient such as
the DMI.

This discussion motivates the following definition.

\begin{definition}
  \label{def:calibration}
  Given $\eta \in C^0(\mathbb{S}^2)$ and
  $\Psi \in C^1(\mathbb{S}^2;\R^2)$ with $\Psi(-e_3) = 0$ and
  $\sup_{z \in \mathbb S^2 \backslash \{-e_3\}} {|\nabla_{\mathbb S^2}
    \Psi(z)| \over |z + e_3|} < \infty$ we call the map
\begin{align}
  \ee_{\eta,\Psi} : \mathcal{A} & \to L^1_{\mathrm{loc}}(\R^2), \notag
  \\
  m  & \mapsto 2 \eta(m) \omega_m + \nabla \cdot \Psi (
       m) \label{eq:lametaPsi} 
\end{align}
the null Lagrangian associated to $\eta$ and $\Psi$.
The corresponding space of null Lagrangians is
\begin{align}
  \begin{split}
    \operatorname{NL} := \Bigg\{ \ee_{\eta, \Psi} : \eta \in
    C^0(\mathbb{S}^2),\, & \Psi\in C^1(\mathbb{S}^2;\R^2) \hbox{ with
    } \Psi(-e_3) = 0 \hbox{ and } \sup_{z \in \mathbb S^2  \backslash \{-e_3\}}
    {|\nabla_{\mathbb S^2} \Psi(z)| \over |z + e_3|} < \infty \Bigg\}.
   \end{split}
\end{align}
\end{definition}

By construction, we have that the integrals of null Lagrangians only
depend on boundary values and the degree.  The growth condition for
$\Psi$ around $-e_3$ allows us to do the standard integration by
parts, noting that $\Psi(m) \in W^{1,1}(\R^2; \R^2)$ for any
$m \in \mathcal A$.
\begin{proposition}
  \label{prop:null_lagrange}
  Let $\eta \in L^\infty(\mathbb{S}^2)$,
  $\Psi \in C^1(\mathbb{S}^2;\R^2)$ with $\Psi(-e_3) = 0$ and
  $\sup_{z \in \mathbb S^2 \backslash \{-e_3\}} {|\nabla_{\mathbb S^2}
    \Psi(z)| \over |z + e_3|} < \infty$.  Then for all
  $m \in \mathcal{A}$ we have
  \begin{align}
    \int_{\R^2} \ee_{\eta, \Psi}(m)
    \intd x  = 2 \mathcal{N}(m) \int_{\mathbb{S}^2} \eta(z)
    \intd \mathcal{H}^2(z). 
  \end{align}
\end{proposition}

A calibration is now simply a null Lagrangian that provides a lower
bound on the energy which is sharp on the presumed minimizer.
\begin{definition}\label{def:calibr}
  We say that a null Lagrangian $\ee \in \operatorname{NL}$ calibrates
  an element $m_* \in \mathcal{A}_1$ (or is a calibration for $m_*$)
  if
\begin{enumerate}[(i)]
\item
  $E_{\bar{\kappa}}(m) \geq \displaystyle \int_{\mathbb{R}^2}
  \lambda(m) \, dx$ \quad for all $m \in \mathcal{A}_1$;
\item
  $E_{\bar{\kappa}}(m_*) = \displaystyle \int_{\mathbb{R}^2}
  \lambda(m_*) \, dx.$
\end{enumerate}
\end{definition}

The main difficulty in proving minimality of a competitor naturally
lies in constructing an appropriate calibration. This is exactly what
we will be doing in the rest of this paper by looking for the
functions $\eta$ and $\Psi$ such that
$e_{\bar\kappa}(m) \geq \lambda(m)$ holds pointwise a.e. in $\R^2$
with $\lambda = \lambda_{\eta,\Psi}$ from
\eqref{eq:lametaPsi}. Needless to say that the right-hand side of the
inequality in Definition \ref{def:calibr} does not depend on $m$ and
is simply equal to $E_{\bar\kappa}(m_*)$, yielding the desired
corollary below.

\begin{corollary}\label{cor:minimality_by_calibration}
  Let $m_* \in \mathcal{A}_1$ be calibrated by
  $\ee \in \operatorname{NL}$. Then $m_*$ is a minimizer of
  $E_{\bar \kappa}$ in $\mathcal{A}_1$.
\end{corollary}

\subsection{Outline of the paper}

We will start by setting up a suitable parametrization of the sphere
$\Sphere^2$ in Section \ref{sec:parametrization}, including a
discussion of equivariant maps in Section \ref{sec:equivariant}.
Next, we construct the null Lagrangian that will turn out to be the
desired calibration in Section \ref{sec:calibration}.  As a warm-up,
we first consider the case of pure exchange energy in Section
\ref{sec:pure_ex}, recasting the Belavin-Polyakov argument in our
language of calibrations.  Section \ref{sec:energy_identity} then
provides the candidate null Lagrangian and the central energy identity
for the full energy.  This is where we carefully systematize and
explain the algebraic manipulations that go into the calibration
construction.  Then Section \ref{sec:calibration_bound} proves that
the candidate null Lagrangian is indeed a calibration via ODE
analysis.  Finally, we conclude in Section \ref{sec:proofs_main} by
proving Theorem \ref{thm:plane}. Unless stated otherwise, proofs of
statements are given at the end of the subsections they appear in.

\section{Parametrizing the sphere and its tangent bundle}\label{sec:parametrization}

We parametrize $ \mathbb{S}^2 \setminus \{ \pm e_3 \} $ via the map
\begin{align}
  z: (0, \pi) \times \mathbb{S}^1 &  \to \mathbb{S}^2 \setminus \{ \pm 
                                    e_3 \}, \notag \\ 
  (\theta, n) & \mapsto - \begin{pmatrix}
    n \sin \theta\\
    \cos \theta
	\end{pmatrix}.
\end{align}
representing the (inverted) spherical coordinates on the sphere. In
particular, $\theta$ is the polar angle measured from the south pole.
For $ (\theta, n) \in (0, \pi) \times \mathbb{S}^1 $, the vectors
corresponding to differentiation along polar and azimuthal direction
on $\mathbb{S}^2$ are given by
\begin{align}
  p := \partial_\theta z
  =
  \begin{pmatrix}
    -n \cos \theta \label{eq:p} \\
    \sin \theta
  \end{pmatrix},\\
  a := -\frac{1}{\sin \theta} \, \partial_n z =
  \begin{pmatrix}
    n^\perp \label{eq:a} \\
    0
  \end{pmatrix},
\end{align}
where $n^\perp := (-n_2,n_1)$ is the counter-clockwise 90$^\circ$
rotation of $n$ and $\partial_n$ is the derivative in the direction
$n^\perp$ on $\mathbb S^1$.  Furthermore, they form an orthonormal
basis of the tangent space $T_{z}\mathbb S^2$ such that the basis
$(z, p, a)$ of $\R^3$ has positive orientation.

Let us record how a Sobolev map $ m \in \mathcal{A}$ behaves when its
target is parametrized in this way.
\begin{lemma}\label{lemma:derivatives_param}
  Let $ m \in \mathcal{A}$.
  For almost all $x\in \R^2$, we define
  \begin{align}
    \theta_m(x) := \arccos (-m_3(x)) \in [0, \pi],
  \end{align}
  and
  \begin{align}
    \label{eq:nm}
    n_m(x) :=
    \begin{cases}
      -\frac{m'(x)}{|m'|(x)} & \text{ if } m(x) \notin \{\pm e_3\},\\
      e_1 & \text{ if }m(x) \in \{\pm e_3\},
    \end{cases}
  \end{align}
  so that almost everywhere, we have the representation
  \begin{align}
    m & = -
	\begin{pmatrix}
          n_m \sin \theta_m\\
          \cos \theta_m
	\end{pmatrix}.
  \end{align}
  
  Furthermore, we have $\theta_m \in H^1(\R^2)$, and with the measure
  theoretic convention $0 \cdot \infty = 0$ the weak chain rule holds:
  \begin{align}
    \nabla \theta_m = \frac{\nabla m_3}{\sin \theta_m}
  \end{align}
  almost everywhere.  Additionally, for $i = 1,2$ we have the
  following identities a. e. in $\R^2$:
  \begin{align}
    \begin{split}
      p \cdot \partial_i m & = \partial_i
      \theta_m, \label{eq:derivative_theta}
    \end{split}\\
    a \cdot \partial_i m & = n_m^\perp \cdot
                           \partial_i  m', \label{eq:aDm} \\
    |\partial_i m|^2 & = |\partial_i \theta_m|^2+ \left|
                       n_m^\perp \cdot \partial_i m'
                       \right|^2. \label{eq:Dm2DthDmp}
  \end{align}
\end{lemma}

We can thus also represent the remaining terms of the energy and the
Jacobian in these coordinates.
\begin{lemma}\label{lemma:energy_param}
  Let $ m \in \mathcal{A}$ and let $\theta_m$ and $n_m$ be as in Lemma
  \ref{lemma:derivatives_param}.  Then we have the following
  representations almost everywhere in $\R^2$:
  \begin{align}
    \omega_m
    & =  \det 
      \begin{pmatrix}
        \partial_1 \theta_m &  \partial_2 \theta_m \\
        n_m^\perp \cdot \partial_1 m' & n_m^\perp \cdot
        \partial_2 m'
      \end{pmatrix}
                                        ,  \label{eq:repr_omega}
    \\ 
    |m'|^2 & = \sin^2 \theta_m,\label{eq:repr_anis}\\
    \begin{split}
      -m' \cdot \nabla m_3
      & 
      = \sin^2  \theta_m (n_m \cdot \nabla \theta_m),
    \end{split}\label{eq:repr_dmi}
  \end{align}
  \begin{align}\label{eq:div_n}
    \nabla \cdot m' = n_m^\perp \cdot  \left( n_m^\perp \cdot
    \nabla m' \right) 
    -\cos \theta_m (n_m \cdot \nabla \theta_m). 
  \end{align}
\end{lemma}

\begin{proof}[Proof of Lemma \ref{lemma:derivatives_param}]
  The representation of $m$ in terms of $\theta_m$ and $n_m$ follows
  straightforwardly from the definitions:
  $|m'|^2 = 1- m_3^2 = 1- \cos^2 \theta_m = \sin^2 \theta_m$ a.e. in
  $\R^2$.  As $\sin \theta_m \geq 0$ almost everywhere, we have
  $|m'| = \sin \theta_m$.
	
  As a result, for $m \in \mathcal A$ we also get
  \begin{align}
    \int_{\left\{\theta_m \leq \frac{\pi}{2}\right\} } \theta^2_m
    \intd x \leq {\pi^2 \over 4} \int_{\left\{\theta_m \leq
    \frac{\pi}{2}\right\}} \sin^2 \theta_m \intd x < \infty. 
  \end{align}
  Furthermore, we have
  \begin{align}
    \int_{\left\{\theta_m > \frac{\pi}{2}\right\} } \theta^2_m  \intd x
    \leq \frac{\pi^2}{4} \int_{\left\{\theta_m  > \frac{\pi}{2}\right\}} \left(1-
    \cos \theta_m  \right)^2  \intd x < \infty, 
  \end{align}
  giving $L^2$-integrability of $\theta_m$.

  For $i = 1,2$ and almost all $x\in \R^2$ such that
  $m(x) \not \in \{\pm e_3\}$, we calculate using
  $m' \cdot \partial_i m' + m_3 \partial_i m_3 = 0$ due to the unit
  length constraint:
  \begin{align}
    \label{eq:1osinDm}
  \begin{split}
    p \cdot \partial_i m & = - n_m \cdot \cos \theta_m \partial_i m' +
    \sin \theta_m \partial_i m_3 = \frac{m' }{|m'|} \cdot \cos
    \theta_m \partial_i m'
    + \sin \theta_m \partial_i m_3\\
    & = - \frac{m_3 \cos \theta_m}{\sin \theta_m} \partial_i m_3 +
    \sin \theta_m \partial_i m_3 = \left(\frac{\cos^2 \theta_m}{\sin
        \theta_m} + \sin \theta_m \right) \partial_i m_3 =
    \frac{1}{\sin \theta_m} \, \partial_i m_3.
  \end{split}
\end{align}
Note that given the convention $0\cdot \infty =0$, this identity also
holds on $\{m = \pm e_3\}$, since $\nabla m = 0$ almost everywhere
thereon.  This also ensures that equation \eqref{eq:aDm} holds almost
everywhere on $\{m = \pm e_3\}$.

For $\eps \in (0, 1)$, let $F_\eps \in C^1([-1,1])$ be defined as
$F_\eps(t) := \arccos((1- \eps)t) \in [0,\pi]$.  Then
$F_\eps(t) \to \arccos t$ uniformly in $t$ as $\eps \to 0$ with
$\sin F_\eps(t) \searrow \sqrt{1 - t^2}$ uniformly in $t$, and
\begin{align}
  F'_\eps(t) = -\frac{1-\eps}{\sqrt{1-(1- \eps)^2t^2}} = -
  \frac{1-\eps}{\sin F_\eps(t)}. 
\end{align}

By the chain rule for Sobolev functions, for $i = 1,2$ we have almost
everywhere:
\begin{align}
  \begin{split}
    \partial_i F_\eps(-m_3) & = \frac{(1-\eps) \partial_i m_3}{\sin
      F_\eps(-m_3)},
  \end{split}
\end{align}
so that by the identity \eqref{eq:1osinDm} and
$\sin F_\eps(-m_3) \geq \sqrt{1-m_3^2} = \sin \theta_m$ we have
\begin{align}
  |\partial_i F_\eps(-m_3)| \leq  \frac{| \partial_i m_3|}{\sin \theta_m} =
  |p \cdot \partial_i m|. 
\end{align}
Since $(z,p,a)$ form an orthonormal basis for all $z\in \Sphere^2$,
for $i = 1,2$ we have a.e.\ on $\{m\neq \pm e_3\}$ that
\begin{align}
  |\partial_i m|^2 = | p \cdot  \partial_i m |^2 +
  \left| 
  n^\perp \cdot \partial_i m' \right|^2.
\end{align}
Note that this also holds on $\{m=\pm e_3\}$ because on this set we
have $\nabla m = 0$ a.e.\ and thus both sides vanish.  Consequently,
we obtain that
\begin{align}
  F_\eps(-m_3) \in  L^2_{\mathrm{loc}}(\R^2) \cap \mathring H^1(\R^2)
\end{align}
is uniformly bounded in $H^1(B_r(0))$ for all $r>0$. In the limit
$\eps \to 0$, we therefore get $F_\eps(-m_3) \to \theta_m$ a.e. and
$\nabla F_\eps(-m_3) \rightharpoonup \nabla \theta_m$ weakly in
$L^2(\R^2)$.  Together with the $L^2$-integrability we established
earlier, we get $\theta_m \in H^1(\R^2)$.
  
To prove the chain rule, note that by $\nabla m_3 = 0$ a.e. on
$\{m=\pm e_3\}$ we get
  \begin{align}
    \nabla F_\eps(- m_3(x)) = \frac{(1-\eps) \nabla m_3(x)}{\sin
    F_\eps(-m_3(x))} \to \frac{\nabla m_3(x)}{\sin \theta_m(x)} 
  \end{align}
  for a.e. $x \in \R^2$ as $\eps \to 0$.  Since we have the uniform,
  $L^2$-integrable dominating function
  \begin{align}
    \left|\nabla F_\eps(- m_3)  - \frac{\nabla m_3}{\sin \theta_m}
    \right| \leq 2\frac{|\nabla m_3|}{\sin \theta_m}, 
  \end{align}
  we get $\nabla F_\eps(- m_3) \to \frac{\nabla m_3}{\sin \theta_m}$
  strongly in $L^2(\R^2)$.  In particular, we get
  $\nabla \theta_m = \frac{\nabla m_3}{\sin \theta_m}$ a.e. from
  uniqueness of the weak limit.
\end{proof}

\begin{proof}[Proof of Lemma \ref{lemma:energy_param}]
  The representation for $\omega_m$ follows from
  \eqref{eq:derivative_theta}, \eqref{eq:aDm} and the fact that
  $(z,p,a)$ is a positively oriented, orthonormal basis of $\R^3$.
  The statements \eqref{eq:repr_anis} and \eqref{eq:repr_dmi}
  immediately follow from Lemma \ref{lemma:derivatives_param}.  On
  $\{m \neq \pm e_3\}$, identity \eqref{eq:div_n} then follows from
  \begin{align}
    \begin{split}
      \nabla\cdot m' & = n_m^\perp \cdot \left( n_m^\perp \cdot \nabla
        m' \right) + n_m \cdot
      \left( n_m \cdot \nabla m'  \right)\\
      & = n_m^\perp \cdot \left( n_m^\perp \cdot \nabla m' \right) +
      \frac{m_3}{|m'|} \left( n_m \cdot
        \nabla m_3 \right) \\
      &= n_m^\perp \cdot \left( n_m^\perp \cdot \nabla m' \right)
      -\cos \theta_m (n_m \cdot \nabla \theta_m),
    \end{split}
  \end{align}
  again, using $m' \cdot \partial_i m' = - m_3 \partial_i m_3$
  for $i = 1,2$. This concludes the proof.
\end{proof}

\section{Radial profiles}\label{sec:equivariant}

As stated in Section \ref{sec:main_results}, radial magnetizations
take the form
\begin{align}
  \m_\theta(x) = -\left(  \sin\theta (|x|) \frac{x}{|x|} , \cos
    \theta(|x|) \right)
\end{align}
for $x\in \R^2$ and a map $\theta: [0, \infty) \to \R$.  For
$\theta \in \mathcal{A}_\mathrm{r}$ (see definition \eqref{eq:Ar}), we
have $m_\theta \in \mathcal{A}$ with
\begin{align}
  n(x) & = \frac{x}{|x|},\\
  \nabla \left( \theta(|x|) \right)& = \theta'(|x|)
                                     n(x), \label{eq:polar}\\ 
  n^\perp(x) \cdot  \partial_i  m'_\theta (x)
       & = - \frac{ \sin 
         \theta(|x|) }{|x|}
         n_i^\perp(x) \qquad i = 1,2, \label{eq:azimuthal} 
\end{align}
for almost all $x\neq 0$.  Consequently, we indeed have
$E_{\bar \kappa}(m_\theta) = I_{\bar \kappa}(\theta)$, where
\begin{align}
  I_{\bar \kappa}(\theta) =2\pi
  \int_0^\infty \left(r |\theta'|^2 + \left(r +\frac{1}{r}
  \right)\sin^2 \theta + 2 \bar \kappa r \theta' \sin^2\theta\right)
  \intd r. 
\end{align}

The degree of $m_\theta$ may be calculated using the following lemma.

\begin{lemma}\label{lemma:degree}
  For $\theta \in \mathcal{A}_\mathrm{r}$, the limits
  $\theta(0):= \lim_{r\to 0}\theta(r) \in \pi \mathbb{Z}$ and
  $\theta(\infty) := \lim_{r\to \infty} \theta(r) \in 2 \pi \mathbb Z$
  exist.  Furthermore, the degree of $m_\theta \in \mathcal A$
  satisfies
	\begin{align}
          \mathcal{N}(m_\theta)= -\frac{1}{2} \int_0^\infty
          \theta'\sin\theta \,
          dr= 
          \frac{1}{2} \big( \cos \theta(\infty) - \cos \theta(0)
          \big)  \in \{0, 1\}.  
	\end{align}
	In particular, we have $m_\theta \in \mathcal{A}_{1}$ if
        $\theta \in \mathcal{A}_{\mathrm{r},1}$.
\end{lemma}

With this setup, we can draw some simple conclusions for a strictly
decreasing critical point $\theta_*$.  First, we state some standard
regularity properties, along with the behavior near $0$ and $\infty$.

\begin{lemma}\label{lemma:reg}
  Let $\bar \kappa > 0 $ support a monotone skyrmion profile
  $\theta_*$ (see
  Definition~\ref{def:monotone_skyr}).  Then we have
  $\theta_* \in C^2([0,\infty))$ with $\theta'_*(0)<0$ and
  $\theta''_*(0) = 0$.  Additionally, there exists $a_*>0$ such that
  in the limit $r\to \infty$ we have
  \begin{align}\label{eq:besseltail}
    \theta_*(r) & = a_* K_1\left(r\right)\left(1+O\left(r^{-3/2} e^{-
                  r}\right) \right),\\ 
    \theta'_*(r) & = a_* K_1'\left(r\right)\left(1+O\left(r^{-3/2}
                   e^{-r}\right) \right), \label{eq:besseltail_deriv} 
  \end{align}
  where $K_1$ is the modified Bessel function of order one.
  Furthermore, as $r\to \infty$ we have
  \begin{align}\label{eq:besselfracfirst}
    \frac{\theta'_*(r)}{\theta_*(r)}
    & = - \left(1+ \frac{1}{2r} +
      O\left(r^{-2}\right)\right) 
  \end{align}
  and $\lim_{r\to \infty} \frac{\theta''_*(r)}{\theta_*(r)} = 1$. In
  particular, we have $\theta_* \in \mathcal{A}_{\mathrm{r},1}$.
\end{lemma}

Next, we express the Euler-Lagrange equation \eqref{eq:ELG} in a way
that will be convenient for what follows.

\begin{lemma}\label{lemma:Pohozhaev}
  Let $\bar \kappa > 0 $ support a monotone skyrmion profile
  $\theta_* \in \mathcal{A}_{\mathrm{r},1}$.  For $r>0$, we introduce
  the quantity
  \begin{align}
    \label{eq:pokh}
    P(r) := \frac{1}{2} \left(\theta'_*(r) \right)^2 - \frac{1}{2} \left(
    1 + \frac{1}{r^2}\right) \sin^2 \theta_*(r).
  \end{align}
  Then $P \in C^\infty(0,\infty)$ and
  \begin{align}\label{eq:derivP}
    r P'(r) = - \left(\theta'_*(r)\right)^2 + \left(
    \frac{1}{r^2} - \bar \kappa  \theta'_* (r) \right) \sin^2
    \theta_*(r), 
  \end{align}
  for all $r>0$.  Furthermore, we have
  \begin{align}
    \lim_{r \to 0} P\left( r \right) = \lim_{r \to \infty} P\left(r\right)  = 0.
  \end{align}
\end{lemma}

Finally, we exploit the fact that $\theta_*$ is strictly monotone
decreasing by introducing its inverse.
\begin{lemma}\label{lemma:rho}
  Let $\bar \kappa > 0$ support a monotone skyrmion profile
  $\theta_* \in \mathcal{A}_{\mathrm{r},1}$.  Then there exists
  $\rho_* \in C^1((0,\pi])$ such that
  \begin{align}
    \theta_* \circ \, \rho_* = \operatorname{id}_{(0,\pi]}, \quad
    \rho_* \circ \, \theta_{*} =  \operatorname{id}_{[0,\infty)}, 
  \end{align}
  and
  \begin{align}
    \label{eq:rhopthp}
    \rho_*'= \frac{1}{\theta'_* ( \rho_*)} < 0.
  \end{align}
\end{lemma}

\begin{proof}[Proof of Lemma \ref{lemma:degree}]
  By Young's inequality, the map $r \mapsto \theta'(r) \sin \theta(r)$
  is integrable over $(0,\infty)$.  Therefore, we have
  $\cos \theta \in W^{1,1}_{\mathrm{loc}}([0, \infty))$.  In
  particular, $a_-:=\lim_{r\to 0} \cos \theta(r)$ exists in $\R$.  Let
  $\delta >0$.  By continuity of $\theta$ on $(0,\infty)$, for $r$
  sufficiently small, $\theta(r)$ cannot escape a connected component
  of $\cos^{-1}(a_- - \delta, a_- + \delta)$.  Therefore, also
  $\theta(0):=\lim_{r\to 0} \theta(r) \in \cos^{-1}(a_-)$.  As a
  result of $\frac{\sin^2 \theta(r)}{r}$ being integrable, we must
  also have
  \begin{align}
    \sin \theta(0) = \lim_{r \to 0} \sin \theta(r) = 0,
  \end{align}
  so that $\theta(0) \in \pi \Z$.
	
  For the limit at infinity, we use the fundamental theorem of
  calculus for Sobolev functions to get
  \begin{align}
    \lim_{r\to \infty} \cos \theta(r) = \lim_{r \to \infty}  \left(
    \cos \theta(1) - \int_1^r \theta'(t)  \sin \theta(t)  dt \right) =
    \cos \theta(1) 
    - \int_1^\infty \theta' \sin \theta  \, dr \in \R,
  \end{align}
  since the last integral converges absolutely.  Integrability of
  $r (1 - \cos \theta)$ yields
  $\lim_{r \to \infty} \cos \theta(r) = 1$.  Therefore, we get
  $\theta(\infty) \in 2 \pi \mathbb Z$, and using the identities
  \eqref{eq:repr_omega}, \eqref{eq:polar}, and \eqref{eq:azimuthal} we
  may calculate
  \begin{align}
    \mathcal{N}(m_\theta) = -\frac{1}{2} \int_0^\infty  \theta'\sin\theta
    dr = 
    \frac{1}{2} \left(\cos \theta(\infty) - \cos \theta(0) \right), 
  \end{align}
  from which $\mathcal{N}(m_\theta) \in \{0, 1\}$ immediately
  follows. Finally, this formula also implies that
  $\mathcal N(m_\theta) = 1$ for
  $\theta \in \mathcal A_{\mathrm r, 1}$, yielding
  $m_\theta \in \mathcal A_1$.
\end{proof}
    
\begin{proof}[Proof of Lemma \ref{lemma:reg}]
    All constants in this proof may depend on $\bar \kappa$ and
    $\theta_*$. 
    
    The regularity statement on $(0,\infty)$ is standard.  Li and
    Melcher \cite[Lemma 14]{li18} provide an argument for the
    asymptotics
    \begin{align}\label{eq:melcher_li_0}
      \theta_*(r) = \pi - b_* r  + o(r),
    \end{align}
    for some $b_* \geq 0$ as $r\to 0$, and
    \begin{align}\label{eq:melcher_li_infinity}
      \theta_*(r) = \frac{c_* e^{-r}}{\sqrt{r}} +
      o\left(\frac{e^{-r}}{\sqrt{r}}  \right) ,
    \end{align}
    for some $c_ * \geq 0$ as $r \to \infty$ (see \cite{dmss:rad} for
    a full treatment).  Consequently, we just need to improve these
    asymptotics to get the desired statement.
	
    \textit{Step 1: Behavior as $r \to 0$.}\\
    Let $\varphi:=\pi-\theta_*$, so we have $\varphi(0)=0$ and
    $\varphi(r) > 0$ for all $r > 0$. As we have
    $\sin\theta_*\cos\theta_*=-\sin\varphi\cos\varphi$ and
    $\sin^2\theta_*=\sin^2\varphi$, for $r \in (0,\infty)$ equation
    \eqref{eq:ELG} becomes
    \begin{align}\label{eq:at_zero_first}
      L_0\varphi =N_0.
    \end{align}
    with
    \begin{align}
      L_0f&:=-(rf')' +  \frac{f}{r},\\
      N_0&:=-\frac{1}{r}\left(\sin\varphi\cos\varphi -
           \varphi\right) - r \sin\varphi\cos\varphi - \bar
           \kappa\sin^2\varphi 
    \end{align}
    for $f \in C^2(0,\infty)$.
				
    The homogeneous equation $L_0f = 0$ on $(0,\eps)$ has the two
    linearly independent solutions, $f(r) = r$ and $f(r) = r^{-1}$. By
    the asymptotics \eqref{eq:melcher_li_0} and variation of
    constants, the equation $L_0 f = N_0$ has the particular solution
    \begin{align}
      f_0 (r) :=- \frac{r}{2}\int_0^r \frac{N_0(t)}{t}\intd
      t+\frac{1}{2r}\int_0^r tN_0(t)\intd t. 
    \end{align}
    One can readily check the bound $|N_0(r)| \leq C r \varphi(r)$ for
    some $C > 0$ universal and all $r>0$ sufficiently small, so that
    from boundedness and monotonicity of $\varphi$ we get
    \begin{align}
      \label{eq:varphi0est}
      \varphi(r) = b_* r + O\left( r^2\varphi(r) \right),
    \end{align}
    as $r \to 0$. In particular, we have $b_* > 0$, since otherwise
    $\varphi(r) \equiv 0$ in some neighborhood of $r = 0$, a
    contradiction to positivity of $\varphi(r)$ for all $r >
    0$. Finally, observe that $f_0(r)$ is twice continuously
    differentiable, with
    \begin{align}
      \label{eq:f0ppN0}
      f_0''(r)=-\frac{N_0(r)}{r}+\frac{1}{r^3}\int_0^rtN_0(t)\,dt,
    \end{align}
    so $\varphi \in C^2([0, \infty))$ and $\theta''(0) = 0$ by
    \eqref{eq:varphi0est}.
						
    \textit{Step 2: Behaviour as $r \to  \infty$.}\\
    We rewrite the Euler-Lagrange equation \eqref{eq:ELG} as
    \begin{align}\label{eq:at_zero_second}
      L_\infty \theta_* =N_\infty,
    \end{align}
    with the modified Bessel operator for $f \in C^2(0,\infty)$ and
    the right-hand side
    \begin{align}
      L_\infty f &:= -f'' - \frac{f'}{r} + \left(1+
                   \frac{1}{r^2}\right) f,\\ 
      N_\infty&:=-\left(1+\frac1{r^2}\right)\left(\sin\theta_*\cos\theta_*- 
                \theta_*\right) 
                + \frac{\bar \kappa}{r}\sin^2\theta_*.
    \end{align}
    The bound \eqref{eq:melcher_li_infinity} gives
    \begin{align}\label{eq:bound_n_infty}
      0 < N_\infty(r) \leq C r^{-2} e^{-2r},
    \end{align}
    for some $C > 0$ and all $r > 0$ sufficiently large.
  
    The homogeneous equation for $L_\infty$ has fundamental solutions
    $I_1$ and $K_1$, the modified Bessel functions of the first and
    second kind of order 1.  The corresponding Wronskian is
 \begin{align}
   K_1(r)  I_1'(r) - I_1(r)K_1'(r)  = {1 \over r},
 \end{align}
 implying that $I_1(r)/K_1(r)$ is increasing. By the asymptotics
 \eqref{eq:bound_n_infty} and variation of constants, a particular
 solution of the equation $L_\infty f = N_\infty$ is given by
 \begin{align}\label{eq:varconst}
   \begin{split}
     f_\infty (r) := I_1\left(r\right)\int_r^\infty t
     K_1\left(t\right)N_\infty(t)\intd t -
     K_1\left(r\right)\int_r^\infty t I_1\left(t\right)N_\infty(t)
     \intd t,
   \end{split}
 \end{align}
 for $r>0$. By an explicit estimate using positivity of $N_\infty(r)$
 and monotonicity of $I_1(r) / K_1(r)$ for $r>0$ large enuogh, the
 first integral is bounded by the second.  Thus, we have
 $f_\infty(r) < 0$ for all large enough $r> 0$. Standard asymptotics
 for Bessel functions together with the asymptotics
 \eqref{eq:melcher_li_infinity} and \eqref{eq:bound_n_infty} give
 \begin{align}
   f_\infty(r) = O\left( r^{-2} e^{-2 r} \right),
 \end{align}
 as $r \to \infty$.  Consequently, we have
 \begin{align}
   \theta_*(r) = a_* K_1(r) + f_\infty(r) = \left\{ a_* +
   O\left(r^{-3/2} e^{-r} 
   \right) \right\} K_1(r) ,
 \end{align}
 for some $a_* \geq 0$ as $r \to \infty$. In particular, we get
 $a_* > 0$, since $f_\infty(r) < 0$ for all $r > 0$ large enough.
 This gives the claim \eqref{eq:besseltail}.
 
 The claim \eqref{eq:besseltail_deriv} follows from the representation
 \eqref{eq:varconst} by differentiating
 $\theta_*(r) = a_* K_1(r) + f_\infty(r)$. Finally, we use the
 Euler-Lagrange equation to get
 \begin{align}
   \frac{\theta''_*}{\theta_*} = - \frac{\theta'_*}{r\theta_*} + 1
   +O\left(r^{-2} + \bar \kappa \theta_* \right), 
 \end{align}
 giving $\lim_{r\to \infty} \frac{\theta''_*}{\theta_*}  = 1$. Lastly,
 regularity and exponential decay of $\theta_*$ give integrability
 and, hence, that $\theta_* \in \mathcal A_{r,1}$. 
\end{proof}

    \begin{proof}[Proof of Lemma \ref{lemma:Pohozhaev}]
      Calculating the derivative of $P$ is a standard computation.  By
      l'H{\^o}pital's rule, we have
      $\lim_{r \to 0} P\left( r \right) = 0$.  The statement
      $ \lim_{r \to \infty} \frac{(\theta'_*)^2(r)}{2} = 0$ is
      contained in Lemma \ref{lemma:reg}.
    \end{proof}
    
    \begin{proof}[Proof of Lemma \ref{lemma:rho}]
    	These are standard statements.
    \end{proof}

\section{Building a calibration for radial
  minimizers}\label{sec:calibration} 

When comparing two maps from $\R^2$ to $\mathbb S^2$ with the aim of
establishing their equality, a natural thing to do is to look at the
difference of their gradients at the same point on the {\em target}
manifold $\mathbb S^2$. These can be written in terms of some basis in
the tangent plane at that point of the target manifold, giving rise to
two vector fields $X, Y : \R^2 \to \R^2$. When one of the maps is the
radial skyrmion, it makes sense to use spherical coordinates and
project onto the directions $p$ and $a$ defined in \eqref{eq:p} and
\eqref{eq:a}, for $x, x_* \in \R^2$ such that $m(x) = m_*(x_*)$, where
$m \in \mathcal A_1$ is a competitor of the radial skyrmion $m_*$.

With the help of the identities \eqref{eq:polar} and
\eqref{eq:azimuthal}, the deviation vector fields $X$ and $Y$ of the
competitor $m \in \mathcal{A}_1$ from the radial skyrmion $m_*$ are
given by
\begin{align}
  X_i & := \partial_i \theta_m - v(\theta_m) n_{m,i}
        , \label{eq:Xi} \\
  Y_i & := n^\perp_m \cdot \partial_i  m' - w(\theta_m)
        n_{m,i}^{\perp}, \label{eq:Yi}
\end{align}
for $i = 1,2$, with two functions $v, w : [0,\pi] \to \R$, provided we
take
\begin{align}
	v(\theta) & := \theta_*' \left( \rho_*(\theta)\right) , \label{eq:v}\\
	w(\theta) & := - \frac{\sin \theta}{\rho_*(\theta)}, \label{eq:w}
\end{align}
extended by continuity to $\theta = 0, \pi$; see Lemma
\ref{lemma:reg_vw}.  Indeed, if we take $m = m_{\theta_*}$, then with
the above choice of $v$ and $w$ we have $X=Y=0$, and this is exactly
what we want from our ``deviation from the radial profile''.

As the difference of energy and calibration should be a coercive
quantity in some sense, let us come up with a quadratic function of
$ X$ and $ Y$ containing information on the original energy and the
topological class. The simple expression
$ \left| X \right|^2 + \left| Y \right|^2 - 2 \eta(m) \det\left( X,
  Y\right)$ satisfies those properties, provided $|\eta(m)|\leq1$.
Pretending to not already know the answer, we consider general
functions $\eta \in L^\infty(\mathbb{S}^2)$ and
$v, w \in L^\infty(0, \pi)$ and almost everywhere on
$\{m \neq \pm e_3\}$ calculate
\begin{align}
  \label{eq:calibgrad} 
    \left| X \right|^2 & + \left| Y \right|^2 - 2 \eta \det\left( X,
      Y\right) \notag \\
    & = |\nabla \theta_m|^2 - 2 v n_m \cdot \nabla \theta_m + v^2 +
    \left| n^\perp_m \cdot \partial_1 m' \right|^2 + \left| n^\perp_m
      \cdot \partial_2 m' \right|^2 - 2 w n_m^\perp \cdot \left(
      n_m^\perp \cdot \nabla m' \right) + w^2
    \notag \\
    &\quad -2 \eta \left[ \det
      \begin{pmatrix}
        \partial_1 \theta_m &  \partial_2
        \theta_m  \\
        n_m^\perp \cdot \partial_1 m' & n_m^\perp \cdot
        \partial_2 m'
      \end{pmatrix}
      + v w - w n_m\cdot \nabla \theta_m - v n_m^\perp \cdot \left(
        n_m^\perp \cdot \nabla m' \right)
    \right] \notag \\
    & = |\nabla \theta_m|^2 + \left| n^\perp_m \cdot
      \partial_1 m' \right|^2 + \left| n^\perp_m \cdot
      \partial_2 m' \right|^2 -2 \eta \det
      \begin{pmatrix}
        \partial_1 \theta_m &  \partial_2 \theta_m \\
        n_m^\perp \cdot \partial_1 m' & n_m^\perp \cdot \partial_2 m'
      \end{pmatrix} \\
      & \quad + 2 \left( v \eta - w \right) n_m^\perp \cdot \left(
          n_m^\perp \cdot \nabla m' \right) +2\left(w \eta -v
      \right) n_m \cdot \nabla \theta_m \notag \\
      & \quad + v^2 + w^2 - 2 vw \eta, \notag 
\end{align}
where we put $\eta = \eta(m)$, $v = v (\theta_m)$ and
$w = w(\theta_m)$ for shorthand.  Here, the first line is the
Dirichlet energy and the Jacobian density in view of Lemmas
\ref{lemma:derivatives_param} and \ref{lemma:energy_param}.  The
second line carries all contributions that are linear in derivatives,
so that these will have to combine to the divergence term in our
calibration, see Definition \ref{def:calibration}.  The last line does
not carry any derivatives, which means that in the whole space
setting, they have to cancel out.  (If we were on a bounded domain,
they would only need to combine to a constant.)

\medskip

The result is the following crucial identity for the gradient of $m$.

\begin{lemma}\label{lemma:energy_id_algebra}
  Let $\eta \in C^0(\mathbb{S}^2)$, $v, w \in C^0([0,\pi])$, and
  $m \in \mathcal{A}$.  Let $X,Y: \R^2 \to \R^2$ be the vector fields
  whose components are defined in \eqref{eq:Xi} and \eqref{eq:Yi}.
  Then
  \begin{align}\label{eq:bp_calibration}
    \begin{split}
      |\nabla m|^2 -2 \eta(m) \omega_m & = \left| X \right|^2 + \left|
        Y \right|^2 -
      2 \eta(m) \det\left( X,  Y\right)\\
      & + 2 \left( w(\theta_m) - v(\theta_m) \eta(m) \right) n_m^\perp
      \cdot \left( n_m^\perp \cdot \nabla m' \right) + 2
      \left(v(\theta_m) - w(\theta_m) \eta(m)
      \right) n_m \cdot \nabla \theta_m \\
      & -\left( v^2(\theta_m) + w^2(\theta_m) -2 v(\theta_m)
        w(\theta_m) \eta(m) \right)
    \end{split}
  \end{align}
  almost everywhere in $\R^2$.
\end{lemma}
We observe that the left-hand side of \eqref{eq:bp_calibration} is
exactly the difference between the energy density and a null
Lagrangian. The right-hand side contains a coercive function (provided
$|\eta(z)| < 1$ for almost all $z \in \mathbb S^2$) plus some terms we
would like to cancel.

\subsection{The case of pure exchange energy}\label{sec:pure_ex}

As a warm-up and to demonstrate the strategy, we want to show how the
idea of calibration naturally appears in the case of pure exchange
energy, meaning we now consider the case of $ \kappa = 0$ and $Q=1$
for the energy in \eqref{eq:E}.  We already mentioned that in this
case Belavin and Polyakov essentially found the calibration
$\ee= 2 \omega_m$, so this is what we will be aiming at now.

Of course, the uniqueness part of Theorem
\ref{thm:existence_mon_radial} does not hold in this case, because
nothing breaks the dilation symmetry of the exchange energy.
Nevertheless, radial critical points exist in the form of the
Belavin-Polyakov profiles given by $m=m_{\theta_*}$ with
\begin{align}\label{def:BP}
	\theta_*(r) =2 \arctan\left( \frac{a}{r}\right),
\end{align}
for $a>0$ and all $r>0$.  In the following, the goal will be to
demonstrate their global minimality in
$\mathring H^1(\R^2;\mathbb{S}^2)$, using the language of
calibrations.  Note the minor notational issue that we do \emph{not}
have $\theta_* \in \mathcal{A}_{r,1}$ due to the usual logarithmic
divergence in the anisotropy energy, which of course plays no role
here.

It is clear that if we choose $v = w$ and $\eta(z) = 1$ for all
$z \in \mathbb{S}^2$, we would make all the spurious terms in identity
\eqref{eq:bp_calibration} vanish. Moreover, with these choices we
would have
\begin{align}
 \begin{split}
   \left| X \right|^2 + \left| Y \right|^2 - 2 \eta(m) \det\left( X,
     Y\right) = \left| X + Y^ \perp \right|^2 \geq 0,
\end{split}
\end{align}
and thus for any $m \in \mathring H^1(\R^2;\mathbb{S}^2)$ with
$\mathcal{N}(m)=1$ we get
\begin{align}\label{eq:energy_identity_heisenberg}
  \begin{split}
    \int_{\R^2} |\nabla m|^2 \intd x - 8\pi & =\int_{\R^2} |\nabla
    m|^2 \intd x - 2 \int_{\R^2} \omega_m \intd x = \int_{\R^2} \left|
      X + Y^ \perp \right|^2 \intd x\geq 0.
  \end{split}
\end{align}
This is precisely the topological bound \eqref{eq:top_bound} for
$\mathcal{N}(m)=1$, as expected.

It therefore remains to see that the identity $v = w$ actually holds
for an admissible profile $\theta_*: [0,\infty) \to [0,\pi]$.
Recalling \eqref{eq:v}, \eqref{eq:w} and \eqref{eq:rhopthp}, the
requirement that $v(\theta)=w(\theta)$ translates into the equation
$ \rho_*'(\theta) =-\rho_*(\theta)/\sin\theta$ for $\rho_*(\theta)$ in
$(0, \pi)$. This equation is uniquely satisfied by
$\rho_*(\theta) = a \cot \frac{\theta}{2}$ with $a > 0$, which is
precisely the inverse of the radial Belavin-Polyakov profile
\eqref{def:BP}.

Thus, we have identified calibrations for the radial Belavin-Polyakov
profiles, giving the bound \eqref{eq:energy_identity_heisenberg}.  It
is also clear that $X + Y^\perp = 0$ with $X$ and $Y$ defined in
\eqref{eq:Xi} and \eqref{eq:Yi} for one Belavin-Polyakov profile also
holds for any other Belavin-Polyakov profile, as all of them have
$\int_{\R^2} |\nabla m|^2 \intd x = 8 \pi$.

\subsection{The energy identity}\label{sec:energy_identity}

Let us now go to more complex situations.  To this end, we consider
$\bar \kappa>0$ which supports a monotone skyrmion profile
$\theta_* \in \mathcal{A}_{\mathrm{r},1}$.

Given $ m \in \mathcal{A}$, define the deviation vector fields $X, Y$
as in \eqref{eq:Xi} and~\eqref{eq:Yi} with the functions $v, w$ from
\eqref{eq:v} and \eqref{eq:w}.  Let us state some of the properties of
$v$ and $w$ from their definitions.

\begin{lemma}\label{lemma:reg_vw}
  Let $\bar \kappa>0$ support a monotone skyrmion profile
  $\theta_* \in \mathcal{A}_{\mathrm{r},1}$ (see
  Definition~\ref{def:monotone_skyr}).  Then we have
  $v,w \in C^1([0,\pi])$. Furthermore, there holds
	\begin{align}
		v(0) & = 0,\quad 
		v(\theta)  < 0,\quad
		v(\pi) = \theta'_*(0),\\
		w(0) & = 0, \quad
		w(\theta) <0,\quad
		w(\pi)  = \theta'_*(0),
	\end{align}
	and
	\begin{align}
          w'(\theta)
          =-\frac{\cos \theta}{\rho_*(\theta)} + \frac{\sin
          \theta}{\rho_*^2(\theta)} \rho_*'(\theta),
	\end{align}
        for all $\theta \in (0, \pi)$.
\end{lemma}

For the full energy and taking advantage of identity \eqref{eq:div_n},
Lemma \ref{lemma:energy_id_algebra} takes the following form:

\begin{lemma}\label{lemma:energy_id_general}
  Let $\bar \kappa>0$ support a monotone skyrmion profile
  $\theta_* \in \mathcal{A}_{\mathrm{r},1}$.  Let
  $m \in \mathcal{A}_1$ and $\eta \in C^0(\mathbb{S}^2)$.  Then
  for almost all $x\in \R^2$ we have
\begin{align}\label{eq:energy_id_general}
  \begin{split}	
    e_{\bar \kappa}(m) - 2 \eta(m)\omega_m & = | X |^2 + |Y
    |^2 - 2 \eta(m) \det(X, Y)\\
    & \quad + 2 \left(w(\theta_m) -v(\theta_m) \eta(m)
    \right) \nabla \cdot m'\\
    & \quad + 2 \left(v(\theta_m) -w(\theta_m) \eta(m)
    \right) (n_m \cdot \nabla \theta_m) \\
    & \quad + 2 \left( \cos \theta_m \left( w(\theta_m) -
        v(\theta_m) \eta(m) \right) + \bar \kappa \sin^2
      \theta_m \right) (n_m \cdot \nabla \theta_m)\\
    & \quad + 2 \left( v(\theta_m) w(\theta_m) \eta(m)
      - P\left( \rho_* (\theta_m) \right) - w^2(\theta_m)
    \right),
  \end{split}
\end{align}
where $e_{\bar\kappa}$ is the energy density defined in
\eqref{eq:energy-density} and $P \in C^0([0, \infty])$ is defined via
\eqref{eq:pokh}.
\end{lemma}

Having dealt with the algebra for general functions
$\eta \in C^0(\mathbb{S}^2)$, it is time to clean up the result by
choosing the correct null Lagrangian.  The strategy to do so is as
follows: As the last line of identity \eqref{eq:energy_id_general}
does not involve any derivatives to be integrated to 0, it must
vanish.  This determines the function $\eta(m) = \eta_*(\theta_m)$,
for a unique choice of $\eta_* = \eta_*(\theta)$.  Furthermore, note
that then the three middle lines above take the form
$2 H(\theta) \nabla \cdot m' + 2 h(\theta) (m'\cdot \nabla \theta)$
for two scalar functions $H$ and $h$.  If using the Euler-Lagrange
equation we can prove that $H' = h$, we have found the required
divergence structure by setting $\Psi(m) = 2 H(\theta_m) m'$.

Before we head into the required algebraic manipulations, let us state
the relevant properties of these auxiliary functions, starting with
$\eta_*$.

\begin{lemma}\label{lemma:eta_simple}
  Let $\bar \kappa>0$ support a monotone skyrmion profile
  $\theta_* \in \mathcal{A}_{\mathrm{r},1}$.  The function
	\begin{align}\label{eq:formula_eta}
          \eta_* (\theta) :=  \frac{
          w^2(\theta) + P( \rho_* (\theta)) }{v(\theta) w(\theta)}   
	\end{align}
	satisfies $\eta_* \in C^1((0,\pi]) \cap C^0([0,\pi])$ with
        $\eta_*(0) = \frac{1}{2}$ and $\eta_*(\pi) = 1$.
\end{lemma}
We also have the following properties of $H$.
\begin{lemma}\label{lemma:endpoint_H}
  Let $\bar \kappa>0$ support a monotone skyrmion profile
  $\theta_* \in \mathcal{A}_{\mathrm{r},1}$.  Let $H :[0,\pi] \to \R$
  be defined as
  \begin{align}\label{eq:formula_H}
    H(\theta) :=w(\theta) - v(\theta) \eta_*(\theta)  =  \frac{ \rho_*
    (\theta)P(\rho_*(\theta))}{\sin \theta}. 
  \end{align}
  Then $H \in C^1([0,\pi])$ with $H(0)=0$.  In particular, the map
  $\Psi: \mathbb S^2 \to \mathbb R^2$ defined by
  \begin{align}\label{eq:formula_Psi}
	\Psi(m) :=2 H(\theta_m) m'
  \end{align}
  is $C^1$ with $\Psi(-e_3) = 0$ and the growth condition
  $|\nabla_{\mathbb S^2} \Psi(z)| \leq C |z +e_3|$ near $z = -e_3$,
  for some $C > 0$ independent of $z \in \mathbb S^2$.
\end{lemma}

With this, we are finally in the position to state the energy identity
for the appropriate choice of null Lagrangian.

\begin{proposition}\label{prop:energy_identity}
  Let $\bar \kappa>0$ support a monotone skyrmion profile
  $\theta_* \in \mathcal{A}_{\mathrm{r},1}$.  Let
  $m \in \mathcal{A}_1$ then for almost all $x \in \R^2$ we have the
  energy equality
  \begin{align}
    \begin{split}
      e_{\bar \kappa}(m) & = | X |^2 + |Y |^2 - 2 \eta_*(\theta_m)
      \det(X, Y) + 2\eta_* (\theta_m) \omega_m + \nabla \cdot \Psi(
      m),
    \end{split}
  \end{align}
  where $\eta_*$ and $\Psi$ are defined in \eqref{eq:formula_eta} and
  \eqref{eq:formula_Psi}, respectively.
\end{proposition}

\begin{proof}[Proof of Lemma \ref{lemma:energy_id_algebra}]
  On $\{ m \not= \pm e_3\}$ the formula holds a. e. by
  \eqref{eq:Dm2DthDmp}, \eqref{eq:repr_omega} and
  \eqref{eq:calibgrad}. At the same time, with our choice $n_m = e_1$
  on $\{ m = \pm e_3\}$ we have $X = - v e_1$, $Y = - w e_2$, and the
  formula holds a. e. in this set too, because $\nabla m = 0$
  a. e. there.
\end{proof}

\begin{proof}[Proof of Lemma \ref{lemma:reg_vw}]
  The fact that $v,w \in L^\infty(0,\pi) \cap C^1(0,\pi)$ follows from
  Lemmas \ref{lemma:reg} and \ref{lemma:rho}.  For $\theta =\pi$, we
  have $\rho_*(\pi) = 0$, so that $v \in C^1((0,\pi])$ also follows
  from Lemma \ref{lemma:reg}.  To handle $w$, we calculate
  \begin{align}
    \lim_{r \to 0} \frac{\sin \theta_*(r)}{r} = -\theta'_*(0).
  \end{align}
  For the derivative, we have
  \begin{align}
    \label{eq:w_deriv_pf_exact}
    w'(\theta)
    = -\frac{\cos \theta}{\rho_*(\theta)} + \frac{\sin
    \theta}{\rho_*^2(\theta)} \rho_*'(\theta), 
  \end{align}
  so that together with $\theta''_*(0) = 0$ we get
  \begin{align}\label{eq:w_deriv_pf}
    w' (\theta_*(r)) =  \frac{1}{r} + \frac{\sin
    \theta_*(r)}{\theta'_*(r)r^2}+ O(r) = o(1),
  \end{align}
  as $r \to 0$.  This ensures that also $w \in C^1((0,\pi])$.
	
  The endpoint $\theta = 0$ corresponds to the limit $r \to \infty$.
  We calculate
  \begin{align}
    v'(\theta) = \theta''_*(\rho_*(\theta) ) \rho_*'(\theta)  =
    \frac{\theta''_*(\rho_*(\theta) )}{\theta'_*(\rho_*(\theta) )},
  \end{align}
  so that Lemma \ref{lemma:reg} gives
  $\lim_{\theta \to 0} v'(\theta) = - 1$, giving $v \in C^1([0,\pi])$.
  From identity \eqref{eq:w_deriv_pf_exact}, we similarly get
  $w\in C^1([0,\pi])$ with $\lim_{\theta \to 0} w'(\theta) = 0$.
\end{proof}

\begin{proof}[Proof of Lemma \ref{lemma:energy_id_general}]
  Inserting the remaining terms of the energy into the result of Lemma
  \ref{lemma:energy_id_algebra} and using Lemma
  \ref{lemma:energy_param}, we get
\begin{align}\label{eq:nomoresignmistake}
  \begin{split}
    & \quad |\nabla \theta_m|^2 + \left| n_m^\perp \cdot \partial_1 m'
    \right|^2 + \left| n_m^\perp \cdot \partial_2 m' \right|^2 + 2
    \bar \kappa \sin^2 \theta_m (n_m \cdot \nabla
    \theta_m) +  \sin^2 \theta_m - 2  \eta(m) \omega_m \\
    & = | X |^2 + |Y |^2 - 2 \eta(m) \det(X, Y)\\
    & \quad + 2 \left(w(\theta_m) - v(\theta_m) \eta(m) \right) \nabla
    \cdot m'\\
    & \quad + 2\left( v(\theta_m)- w(\theta_m) \eta(m) + \cos \theta_m
      \left(w(\theta_m) - v(\theta_m) \eta(m) \right)+ \bar \kappa
      \sin^2 \theta_m \right)(n_m \cdot \nabla
    \theta_m) \\
    & \quad + \sin^2 \theta_m - \left(v^2(\theta_m) + w^2(\theta_m) -
      2 v(\theta_m) w(\theta_m) \eta(m) \right).
  \end{split}
\end{align}
For the last line, we additionally recall the definition of $P$ in
\eqref{eq:pokh}, together with \eqref{eq:v} and \eqref{eq:w}, and
calculate
\begin{align}
  \begin{split}
    \sin^2 \theta_m - \left(v^2(\theta_m) + w^2(\theta_m) -
      2 v(\theta_m) w(\theta_m) \eta(m) \right) = 2
    \left( v(\theta_m) w(\theta_m) \eta(m) - P\left(
        \rho_* \right) - w^2(\theta_m) \right),
  \end{split}
\end{align}
concluding the proof.
\end{proof}

\begin{proof}[Proof of Lemma \ref{lemma:eta_simple}]
  By Lemmas \ref{lemma:reg}, \ref{lemma:Pohozhaev}, \ref{lemma:rho},
  and \ref{lemma:reg_vw}, we have $\eta_* \in C^1((0,\pi])$ with
  $\eta_*(\pi) = 1$.  In the limit $\theta \to 0$, we however have
  $w(0) = 0$, so that we have to investigate further.
	
  Note that by Lemma \ref{lemma:reg}
	\begin{align}
          \lim_{\theta \to 0} \frac{w}{v} = - \lim_{r\to \infty} 
          \frac{ \sin \theta_*(r)}{r \theta'_*(r)}= 0,
	\end{align}
        and hence
	\begin{align}\label{eq:reuse_eta}
	  \begin{split}
            \lim_{\theta \to 0} \eta_*(\theta) & =- \lim_{r \to
              \infty} \frac{rP(r)}{\theta'_*(r) \sin
              \theta_*(r)}\\
            & = - \lim_{r \to \infty} \left(
              \frac{r\theta'_*(r)}{2\sin \theta_*(r)} -
              \left( \frac{1}{r} + r\right) \frac{\sin
                \theta_*(r)}{2\theta'_*(r)} \right)\\
            & = - \lim_{r \to \infty}\frac{r\theta'_*(r)}{2 \sin
              \theta_*(r)} \left( 1 -\frac{ \theta_*^2
                (r)}{(\theta'_*(r))^2} \right).
	  \end{split}
	\end{align}
	Identity \eqref{eq:besselfracfirst} thus gives
	\begin{align}\label{eq:limit_eta_infty}
	  \begin{split}
		\lim_{\theta \to 0} \eta_*(\theta) & = \frac{1}{2},
	  \end{split}
	\end{align}
	concluding the proof.
\end{proof}

\begin{proof}[Proof of Lemma \ref{lemma:endpoint_H}]
  For $\theta \in (0,\pi)$, we have
  \begin{align}
    H(\theta) =w(\theta) - v(\theta) \eta_{*}(\theta) =
    - \frac{ P(\rho_*(\theta)) }{w(\theta)} = 
    \frac{ \rho_*(\theta) P(\rho_*(\theta)) }{\sin \theta},  
  \end{align}
  and we know that $H \in C^1((0,\pi])\cap C^0([0,\pi])$ and
  $H(0) = 0$ by Lemmas \ref{lemma:reg_vw} and \ref{lemma:eta_simple}.
	
  To handle $H'(0)$, we calculate
  \begin{align}
    H'(\theta) & =  \frac{\rho_*'(\theta)  P(\rho_*(\theta))
                 }{\sin \theta} + 
                 \frac{ \rho_*(\theta) \rho_*'(\theta)
                 P'(\rho_*(\theta))}{\sin \theta} - 
                 \frac{\rho_*(\theta)
                 P(\rho_*(\theta)) \cos \theta }{\sin^2 \theta}.  
  \end{align}
  Hence by Lemma \ref{lemma:rho}
  \begin{align}
    \label{eq:Hp0}
    \lim_{\theta \to 0} H'(\theta)
    & =  \lim_{r \to \infty} \left( \frac{
      P(r) }{\theta_*'(r) \sin \theta_*(r)} + 
      \frac{ r
      P'(r)}{\theta_*'(r) \sin \theta_*(r)} - 
      \frac{r
      P(r) \cos \theta_*(r) }{\sin^2 \theta_*(r)} \right).  
  \end{align}
  By \eqref{eq:reuse_eta} and \eqref{eq:limit_eta_infty}, the first
  term drops out in the limit $r\to \infty$.  For the remaining terms,
  note that we still have
  \begin{align}
    \sin \theta_*(r) & =  \theta_*(r) + O(\theta_*^3(r)),\\
    \cos \theta_*(r) & = 1 + O(\theta_*^2(r)).
  \end{align}
  As $\theta_*(r)$, $P(r)$, and $P'(r)$ decay exponentially fast as
  $r\to \infty$, these expansions give
  \begin{align}
    \begin{split}
      \lim_{\theta \to 0} H'(\theta) & = \lim_{r \to \infty}
      \left(\frac{ r P'(r)}{ \theta'_*(r) \theta_*(r)} -
        \frac{ r P(r)}{ \theta_*^2 (r)} \right).
    \end{split}
  \end{align}
  We now insert the expressions from Lemma \ref{lemma:Pohozhaev} and
  use the expansions of Lemma \ref{lemma:reg} to get
  \begin{align}
    \begin{split}
      \lim_{\theta \to 0} H'(\theta) & = \lim_{r \to \infty} \left[
        -\frac{\theta'_*(r) }{ \theta_*(r)}+
        \frac{r}{2}\left( 1 - \frac{ (\theta'_*(r))^2}{
            \theta_*^2 (r)} \right)\right] =
      \frac{1}{2}.
    \end{split}
  \end{align}

  To conclude the proof, we note that by definition $\Psi(-e_3) = 0$
  and with $\theta_z = \arccos(-z_3)$ we have
  $|\nabla_{\mathbb S^2} \Psi(z)| \leq C (|H(\theta_z)| +
  |H'(\theta_z)| \sin \theta_z) \leq C' |z'| \leq C' | z + e_3|$ with
  some $C, C' > 0$ and all $z \in \mathbb S^2$ near $-e_3$, as
  required.
\end{proof}	

\begin{proof}[Proof of Proposition \ref{prop:energy_identity}]
  Using the formula \eqref{eq:derivP} for the derivative of $P$, on
  $(0,\pi)$ we get
  \begin{align}
    \begin{split}
      v(\theta) - w(\theta) \eta_{*}(\theta) & =
      \frac{v^2(\theta) - w^2(\theta) -
        P(\rho_*(\theta))  }{v(\theta)} \\
      & =- \frac{\rho_*(\theta) P'(\rho_*(\theta)) + \bar
        \kappa \theta'_*(\rho_*(\theta)) \sin^2 \theta +
        P(\rho_*(\theta))
      }{\theta'_*(\rho_*(\theta))}\\
      & = - \left(\rho_*(\theta) P(\rho_*(\theta))\right)' -
      \bar \kappa \sin^2\theta.
    \end{split}
  \end{align}
  Therefore, recalling identity \eqref{eq:formula_H}, we get
  \begin{align}
    \begin{split}
      v(\theta)- w(\theta) \eta_*(\theta) + & \cos
      \theta \left(w(\theta) - v(\theta) \eta_*(\theta)
      \right)+ \bar \kappa \sin^2 \theta \\
      & = - \left(\rho_*(\theta) P(\rho_*(\theta))\right)' +
      \frac{ \cos \theta }{\sin \theta}\rho_*(\theta)
      P(\rho_*(\theta)) \\
      & = - \sin \theta \left( {\rho_*(\theta)
          P(\rho_*(\theta)) \over \sin \theta} \right)'\\
      & = - H'(\theta) \sin \theta .
    \end{split}
  \end{align}
  Inserting everything into the expression
  \eqref{eq:energy_id_general}, we get
  \begin{align}
    \begin{split}
      e_{\bar \kappa}(m) - 2 \eta_*(\theta_m) \omega_m & = | X
      |^2 + |Y |^2 - 2 \eta_*(\theta_m) \det(X, Y) \\
      & \quad + 2 H (\theta_m) \nabla \cdot m' + 2 H' (\theta_m) (m'
      \cdot \nabla \theta_m)
    \end{split}
  \end{align}
  almost everywhere, which yields the statement.
\end{proof}
	
\section{The calibration bound
  $\|\eta_{*}\|_\infty \leq 1$}\label{sec:calibration_bound}

\begin{figure}[t]
  \centering
  \includegraphics[width=8cm]{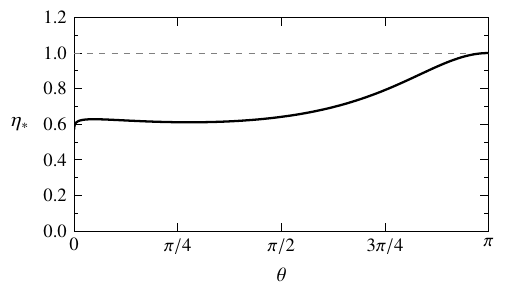}
  \caption{The function $\eta_*(\theta)$ obtained from the numerical
    skyrmion profile for $\bar\kappa = 1$.}
  \label{fig:eta}
\end{figure}

In Section \ref{sec:equivariant}, we have seen that the functions $v$
and $w$ parametrize $\nabla m_*$, and in Section
\ref{sec:energy_identity}, we have uncovered their structural
connection to the exchange energy of any competitor.  However, we have
yet to prove that our lower bound is coercive.  In other words, we
need the crucial calibration bound $\|\eta_{*}\|_\infty \leq 1$ to
show the energy bound in the Definition~\ref{def:calibr} and conclude
uniqueness. We note that a quick way to convince oneself is to plot
the function $\eta_*(\theta)$, using a numerical radial skyrmion
profile, which can be readily computed from \eqref{eq:ELG} to high
accuracy by the shooting method. An example of such a plot for
$\bar \kappa = 1$ is shown in Fig. \ref{fig:eta}.

We recall that
\begin{align}
  \eta_* (\theta) &=  \frac{w^2(\theta) + P(\rho_* (\theta))}{v
                    (\theta) w(\theta)} ,\\  
  v(\theta) & = \theta_*' \left( \rho_*(\theta)\right),\\
  w(\theta) & = - \frac{\sin \theta}{\rho_*(\theta)}
\end{align}
where
\begin{align}
  P(\rho_*(\theta)) = \frac12 \left(\theta'_* (\rho_*(\theta)) \right)^2
  - \frac12 \left(1 + \frac{1}{\rho_*^2(\theta)} \right) \sin^2 \theta, 
\end{align}
and we have used $\theta_*(\rho_*(\theta))=\theta$.

Expressing $\eta_*$ purely in terms of $v$ and $w$ gives
\begin{align}
  \eta_*(\theta) = \frac{1}{2}\left(\frac{v(\theta)}{w(\theta)} +
  \frac{w(\theta)}{v(\theta)} - \frac{\sin^2
  \theta}{v(\theta)w(\theta)} \right), 
\end{align}
which looks like it has a special structure.  In particular, we see
that the calibration bound $|\eta_{*}| \leq 1$ allows us to complete
the squares to equivalently be
\begin{align}\label{eq:calib_bogo}
  \left(\sqrt{\frac{v(\theta)}{w(\theta)}} -
  \sqrt{\frac{w(\theta)}{v(\theta)}} \right)^2 \leq 
  \frac{\sin^2 \theta}{v(\theta) w(\theta)} \leq
  \left(\sqrt{\frac{v(\theta)}{w(\theta)}} + 
  \sqrt{\frac{w(\theta)}{v(\theta)}} \right)^2.  
\end{align}

We therefore have to estimate the two homogeneous quantities of two
variables outside against the quantity in the middle.  This is
naturally done on the quotient of the variables. Because having
$v(\theta)=\theta'_*(\rho_*(\theta))$ in the numerator seems nicer, we
pick the positive quantity $\frac{v(\theta)}{w(\theta)}$ and multiply
it with the above inequality. Introducing the change of variables
$r=\rho_*(\theta)$ and evaluating the resulting expression on
$\theta = \theta_*(r)$ we obtain
\begin{align} \label{eq:cal} \left(\frac{v (\theta_*(r))}{w
      (\theta_*(r))} -1 \right)^2 \leq \frac{\sin^2
    \theta_*(r)}{w^2(\theta_*(r))} \leq \left(\frac{v (\theta_*(r))}{w
      (\theta_*(r))} + 1\right)^2,
\end{align}
where the expression in the middle is simply equal to $r^2$.  In terms
of the scale-free quantity
\begin{align}
  q(r) :=  \frac{v (\theta_*(r))}{w (\theta_*(r))} =
  -\frac{r\theta_*'(r)}{\sin\theta_*(r)} 
\end{align}
for $r \in (0, \infty)$, we can therefore rewrite the inequality \eqref{eq:cal} as
\begin{align}\label{eq:scale-free_top}
	\left(q(r) -1 \right)^2 \leq r^2  \leq \left(q(r)+ 1\right)^2.
\end{align}

In order to interpret this inequality, let us observe that Section
\ref{sec:energy_identity} identifies $v$ and $w$ as good variables for
the skyrmion problem \emph{when evaluated on level sets}.  This
suggests rewriting the radial energy $I_{\bar \kappa}(\theta_*)$ as
\begin{align}
 \begin{split}
   I_{\bar \kappa}(\theta_*) & = 2\pi \int_0^\infty \left(r
     |\theta'_*|^2 + \left(r +\frac{1}{r} \right)\sin^2 \theta_* + 2
     \bar \kappa r \theta'_* \sin^2\theta_*\right) \intd r\\ 
   & = 2\pi \int_{0}^\pi \left(\rho_* |v| + \left(\rho_*
       +\frac{1}{\rho_*} \right)\frac{\sin^2 \theta}{|v|} - 2 \bar
     \kappa \rho_* \sin^2\theta\right) \intd \theta\\ 
   & = 2\pi \int_{0}^\pi\sin \theta \left(\frac{v}{w} + \frac{\sin^2
       \theta}{vw} + \frac{w}{v} + 2 \bar \kappa \frac{\sin^2
       \theta}{w}\right) \intd \theta.
 \end{split}
\end{align}
One then recognizes $\frac{v}{w} + \frac{w}{v}$ as the Dirichlet
energy and $\frac{\sin^2 \theta}{vw}$ as the anisotropy term, both
divided by $\sin \theta$.  Therefore, we can reformulate the strict
calibration bounds
\begin{align}
  -1< \frac{1}{2}\left(\frac{v(\theta)}{w(\theta)} +
  \frac{w(\theta)}{v(\theta)} - 
  \frac{\sin^2 \theta}{v(\theta) w(\theta)} \right) < 1, 
\end{align}
for $\theta \in (0, \pi)$ by undoing the above procedure and
multiplying by $-v(\theta) \sin \theta > 0$ as
\begin{align}
  \theta'_*(r) \sin \theta_*(r)  < \frac{1}{2}\left(r
  |\theta'_*(r) |^2 + 
  \frac{\sin^2\theta_*(r) }{r} - r \sin^2 \theta_*(r)
  \right) < -\theta'_*(r)  
  \sin \theta_*(r),
\end{align}
for all $r > 0$ by undoing the above procedure and multiplying by
$-v(\theta) \sin \theta > 0$.

In view of Lemma \ref{lemma:degree}, the outer functions in this
inequality are simply $\mp r \omega_{m_*}$.  Even better, this means
that completing the square in the inequality \eqref{eq:calib_bogo} is
\emph{precisely} the trick Belavin and Polyakov used to prove the
topological bound \eqref{eq:top_bound}. Then the first inequality in
\eqref{eq:scale-free_top} is a scale-free way to express the defect in
the topological bound \eqref{eq:top_bound} with the choice
$+8\pi \mathcal{N}(m)$, while the second inquality expresses the
defect for $-8\pi \mathcal{N}(m)$.  In particular, the first
inequality controls how far $m$ is from being a Belavin-Polyakov
profile, while the second ensures that $m$ is far away from an
``anti''-Belavin-Polyakov profile.  Therefore, our proof has to rely
on the chiral nature of the DMI term.

On a technical level, the inequality \eqref{eq:scale-free_top} amounts
to comparing $q$ with the three affine functions $1+r$, $1-r$ and
$r-1$. Each of the three differences satisfies a linear ODE of the
form $h'=a\,h+\sin\theta_*\,F(\theta_*)$, see \eqref{eq:gap_crossing},
in which the inhomogeneity $F$ is strictly increasing in $\theta$ and
hence strictly decreasing in $r$.  The inequality then follows from
basic comparison arguments.  This only uses the Euler--Lagrange
equation \eqref{eq:ELG}, strict monotonicity of $\theta_*$, and the
endpoint estimates of Lemma~\ref{lemma:reg}.

\begin{proposition}\label{prop:calib_bound_R2}
Let $\bar\kappa>0$ support a monotone skyrmion
$\theta_*\in\mathcal A_{\mathrm r,1}$. Then
\begin{equation}\label{eq:scale-free_top2}
  |r-1|<q(r) <r+1
  \qquad\text{for all }r>0.
\end{equation}
Consequently, $|\eta_*(\theta)|\leq1$ for all
$\theta\in[0,\pi]$, with equality only at $\theta=\pi$.
\end{proposition}

\begin{proof}
  Since $\theta_*'(r) <0$ and $0<\theta_*(r) <\pi$, we have
  $q>0$.  We have seen above that $|\eta_*(\theta)|<1$ is equivalent
  to $(q(r) -1)^2<r^2<(q(r) +1)^2$, or to
  \eqref{eq:scale-free_top2}.  By Lemma~\ref{lemma:reg}, we have
\begin{equation}\label{eq:q_endpoints}
  q(r)=1+o(r)\quad \text{ as } r \to 0,
  \qquad
  q(r)=r+\tfrac12+O(r^{-1})\quad \text{ as } r\to\infty.
\end{equation}
Differentiating $q$ and using \eqref{eq:ELG} yields
\begin{equation}\label{eq:riccati}
  q'(r) =\frac{\cos\theta_*(r) }{r}(q^2(r) -1-r^2)
  +\bar\kappa\sin\theta_*(r),
\end{equation}
which in particular gives $q'(0)  = 0$.

We now define three gap functions $h_1(r):= 1+r-q(r)$,
$h_2(r):=q(r)+r-1$, $h_3(r):=q(r)+1-r$, and want to show that all
these functions are non-negative.  We will prove afterwards that this
ensures \eqref{eq:scale-free_top2}.  We note that by
\eqref{eq:q_endpoints} and \eqref{eq:riccati} all three gaps are
positive near $0$ and $\infty$: For $i =1,2$, we have $h_i(0) = 0$ and
$h'_i(0) = 1$, so that $h_i(r)>0$ for $r>0$ small enough.  The other
statements follow similarly.

To extend these statements to all $r>0$, we observe that for each gap
function $h_i(r)$ the identity \eqref{eq:riccati} gives
\begin{equation}\label{eq:gap_crossing}
  h_i'(r) =a_i(r)h_i(r) +\sin\theta_*(r)
  \,F_i(\theta_*(r) ), 
\end{equation}
where $a_i(r)$ and $F_i(\theta_*(r))$ are defined
as:
\[
\renewcommand{\arraystretch}{1.8}
\begin{array}{c|c|c|c}
  i& h_i(r)  & a_i(r) & F_i(\theta_*(r) ) \\ \hline
  1&  1+r-q(r)  & \dfrac{\cos\theta_*(r) }{r}(q(r) +1+r)
                           &
                             \dfrac{1-2\cos\theta_*(r)
                             }{\sin\theta_*(r) }-\bar\kappa \\  
  2& q(r) +r-1 & \dfrac{\cos\theta_*(r) }{r}(q(r) +1-r)
                           & \dfrac{1-2\cos\theta_*(r)
                             }{\sin\theta_*(r) }+\bar\kappa \\ 
  3&  q(r) +1-r & \dfrac{\cos\theta_*(r) }{r}(q(r) -1+r)
                           & \dfrac{-1-2\cos\theta_*(r)
                             }{\sin\theta_*(r) }+\bar\kappa 
\end{array}
\]
Moreover, each $F_i(\theta)$ is strictly
increasing in $\theta$, since
\begin{align}
  \frac{d}{d\theta} \left( 
  \frac{t-2\cos\theta}{\sin\theta} \right)
  =\frac{2-t\cos\theta}{\sin^2\theta}>0,
  \qquad t \in\{-1,1\}.
\end{align}
As we know that $\theta_*(r)$ is strictly decreasing we have
$r\mapsto F_i(\theta_*(r))$ is strictly decreasing.

Suppose a gap function $h_i$ is negative on a connected component
$(\alpha,\beta)$ of the set $\{r \in (0,\infty): h_i(r)<0\}$. We
observe that since $h_i>0$ near $0$ and $\infty$, there exists a
component with $0<\alpha<\beta<\infty$, with
$h_i(\alpha)=h_i(\beta)=0$ and $h_i'(\alpha)\leq0\leq
h_i'(\beta)$. Equation~\eqref{eq:gap_crossing} then implies
\[
  F_i(\theta_*(\alpha))\leq0\leq F_i(\theta_*(\beta)),
\]
contradicting strict decrease. Hence $h_i\geq0$ everywhere.  At any
interior zero, $h_i(r_0)=h_i'(r_0)=F_i(\theta_*(r_0))=0$, so
differentiating \eqref{eq:gap_crossing} gives
\begin{align}
  h_i''(r_0)=\sin\theta_*(r_0)\,F_i'(\theta_*(r_0))\theta_*'(r_0)<0,   
\end{align}
contradicting a local minimum condition. All three gaps are therefore strictly
positive, proving \eqref{eq:scale-free_top2}. Finally,
$\eta_*(0)=\tfrac12$ and $\eta_*(\pi)=1$ by
Lemma~\ref{lemma:eta_simple}, completing the proof.
\end{proof}

\section{Proof of the main theorem}\label{sec:proofs_main}

We start by ruling out existence of minimizers for the critical value
of $\bar\kappa = 4 / \pi$.

  \begin{proposition}
    \label{prop:crit}
  For $\bar \kappa = 4/\pi$, we have
  $\inf_{\mathcal{A}_1} E_{4/\pi} = 0$ and $E_{4/\pi}(m) >0$ for all
  $m \in \mathcal{A}_1$.  In particular, the infimum is not attained.
\end{proposition}

\begin{proof}
  Let $m\in\mathcal A_1$. We use the choice $v(\theta)=-\sin\theta$
  and $w(\theta)=0$ in Lemma~\ref{lemma:energy_id_algebra}, motivated
  by asymptotic skyrmion profile as $\bar \kappa \to 4/\pi$
  \cite{dmss:rad}. The corresponding deviation fields are
\begin{align}
  X=\nabla\theta_m-m',\qquad
  Y_i=n_m^\perp\cdot\partial_i m',\quad i=1,2.
\end{align}
As in \eqref{eq:nomoresignmistake}, with
$\eta(m)=\eta_*(\theta_m)$, we obtain
\begin{align}
  \begin{split}
    &e_{4/\pi}(m)-2\eta(m)\omega_m
    -|X|^2-|Y|^2+2\eta(m)\det(X,Y) \\
    &\qquad=2\eta(m)\sin\theta_m\,\nabla\cdot m'\\
    &\qquad\quad+2\left(1-\eta(m)\cos\theta_m
      -\frac4\pi\sin\theta_m\right)
    m'\cdot\nabla\theta_m\\
    &\qquad=\nabla\cdot\Psi(m),
  \end{split}
\end{align}
with $\Psi(m)=2\eta_*(\theta_m)\sin\theta_m\,m'$, provided that
$\eta_*$ is chosen to satisfy
\begin{align}
  \eta_*'(\theta)\sin\theta+2\eta_*(\theta)\cos\theta
  =1-\frac4\pi\sin\theta.
\end{align}
Such a choice is explicitly
\begin{align}
  \eta_*(\theta):=
  \frac{1-\cos\theta-\frac2\pi\theta+\frac1\pi\sin 2\theta}
       {\sin^2\theta},
  \quad 0<\theta<\pi.
\end{align}
In particular, the endpoint singularities are removable, and $\eta_*$
extends smoothly to $[0,\pi]$, with
\begin{align}
  \eta_*(0)=\frac12,\qquad \eta_*(\pi)=-\frac12,
  \qquad \eta_*(\pi-\theta)=-\eta_*(\theta).
\end{align}

We claim that $|\eta_*|\leq\frac12$ on $[0,\pi]$.  Indeed,
\begin{align}
  (\eta_*(\theta)\sin^2\theta)'
  =\sin\theta\left(1-\frac4\pi\sin\theta\right)
\end{align}
changes sign exactly once on $(0,\pi/2)$, from positive to negative.
Since $\eta_*(\theta)\sin^2\theta$ vanishes at both endpoints of
$[0,\pi/2]$, it follows that $\eta_*\geq0$ on this interval. Also,
\begin{align}
  \left[\left(\frac12-\eta_*(\theta)\right)\sin^2\theta\right]'
  =\sin\theta\left(\cos\theta-1+\frac4\pi\sin\theta\right)
  \geq0,
\end{align}
since $4/\pi>1$ and $\sin\theta+\cos\theta\geq1$ on this interval.
The expression in brackets vanishes at $\theta=0$, so
$\eta_*\leq\frac12$ on $[0,\pi/2]$.
The claim follows by antisymmetry.

Writing $\theta_z:=\arccos(-z_3)$ for $z\in\mathbb S^2$, we have
$\eta(z)=\eta_*(\theta_z)\in C^0(\mathbb S^2)$, and by symmetry
\begin{align}
  \int_{\mathbb S^2}\eta(z)\intd\mathcal H^2(z)
  =2\pi\int_0^\pi\eta_*(\theta)\sin\theta\intd\theta=0.
\end{align}
The map $\Psi(z)=2\eta_*(\theta_z)\sin\theta_z\,z'$ extends to a $C^1$
map across both poles, with $\Psi(-e_3)=0$, and satisfies
\begin{align}
  |\nabla_{\mathbb S^2} \Psi(z)|
  \leq C \left( |\eta_*(\theta_z)|\sin\theta_z
  +|\eta_*'(\theta_z)|\sin^2\theta_z \right)
  \leq C'|z+e_3|,
\end{align}
for some $C,C' > 0$ universal.  Thus $\eta$ and $\Psi$ satisfy the
hypotheses of Proposition~\ref{prop:null_lagrange}. Integrating the
energy identity above and using that proposition, together with
$2|\det(X,Y)|\leq|X|^2+|Y|^2$, we conclude that
\begin{align}
  E_{4/\pi}(m)
  \geq\frac12\int_{\R^2}\left(|X|^2+|Y|^2\right)\intd x
  \geq0.
\end{align}
If equality held, then $Y=0$ almost everywhere. By
\eqref{eq:repr_omega}, this would imply $\omega_m=0$ almost
everywhere, contradicting $\mathcal N(m)=1$. Hence
$E_{4/\pi}(m)>0$ for every $m\in\mathcal A_1$.

Finally, the radial result in \cite{dmss:rad} supplies a sequence
$\theta_k\in\mathcal A_{\mathrm r,1}$ with
$I_{4/\pi}(\theta_k)\to0$. Since $m_{\theta_k}\in\mathcal A_1$ and
$E_{4/\pi}(m_{\theta_k})=I_{4/\pi}(\theta_k)$, this gives the matching
upper bound. Therefore $\inf_{\mathcal A_1}E_{4/\pi}=0$, and the
infimum is not attained.
\end{proof}

We finally prove Theorem \ref{thm:plane}.  The strategy is to use the
coercivity properties of the energy identity, Proposition
\ref{prop:energy_identity}, to prove that the map
$\Phi(x):= x - \rho_*(\theta_m(x)) n_m(x)$ is constant if
$m \in \mathcal{A}_1$ is a minimizer, so that $\rho_*(\theta_m)$ is
the distance to a unique point.

\begin{proof}[Proof of Theorem \ref{thm:plane}]
  Smoothness of $\theta_*$ is contained in Lemma \ref{lemma:reg}.
  Combining Corollary \ref{cor:minimality_by_calibration}, the energy
  identity from Proposition \ref{prop:energy_identity}, the
  calibration bound from Proposition \ref{prop:calib_bound_R2}, and
  the bound \eqref{eq:jac-bound} for the Jacobian in terms of
  $|\nabla m|^2$, we see that $m_*$ is a minimizer of
  $E_{\bar \kappa}$ over $\mathcal{A}_1$.
	
  Let $m \in \mathcal{A}_1$ be another minimizer.  By the regularity
  theory \cite{coifman93,helein,chang99} for two-dimensional problems
  of harmonic-map type it is smooth.  In particular, defining the
  non-empty, open set $\Omega_m:= \{m \neq \pm e_3\}$, the maps
  $\theta_m$ and $n_m$ are smooth on $\Omega_m$.  As the calibration
  bound is only saturated at $\theta = \pi$, we have
  $ \nabla \theta_m = v(\theta_m) n_m$ and
  $n^\perp_m \cdot \partial_i m' = w(\theta_m) n_{m,i}^{\perp}$,
  $i = 1,2$, everywhere on $\Omega_m$.
	
  Inserting $ \nabla \theta_m = v(\theta_m) n_m$ and the definition of
  $v$, we calculate
  \begin{align}
    \nabla \rho_*(\theta_m) = \rho_*'(\theta_m) \nabla \theta_m
    =  \rho_*'(\theta_m) v(\theta_m) n_m = n_m 
  \end{align}
  on $\Omega_m$.  Additionally, we can insert
  $n_m^\perp \cdot \partial_i m' = w(\theta_m)n_{m,i}^\perp$ on
  $\Omega_m$, for $i=1,2$, to get
\begin{align}
  \partial_i n_m
  &= \left(n_m^\perp \cdot \partial_i n_m\right)n_m^\perp \\
  &= -\frac{1}{|m'|}
     \left(n_m^\perp \cdot \partial_i m'\right)n_m^\perp \\
  &= -\frac{w(\theta_m)}{|m'|}
     n_{m,i}^\perp n_m^\perp.
\end{align}
On $\Omega_m$, the definition of $w$ therefore gives
\begin{align}
  \partial_i n_m
  = \frac{1}{\rho_*(\theta_m)}n_{m,i}^\perp n_m^\perp.
\end{align}

Considering $\Phi(x):=x-\rho_*(\theta_m(x))n_m(x)$ for
$x\in\Omega_m$, and using
$\partial_i[\rho_*(\theta_m)]=n_{m,i}$, we obtain
\begin{align}
  \partial_i\Phi
  &= e_i-\partial_i[\rho_*(\theta_m)]n_m
       -\rho_*(\theta_m)\partial_i n_m \\
  &= e_i-n_{m,i}n_m-n_{m,i}^\perp n_m^\perp
   =0,
  \qquad i=1,2.
\end{align}
Thus $\nabla\Phi=0$ on $\Omega_m$.
  On each connected component $U$ of $\Omega_m$ we therefore have
  $\Phi = x_U$ for some $x_U \in \R^2$.  For every $x\in U$, we have
  $x - x_U = \rho_*(\theta_m) n_m$, so that
  $\rho_*(\theta_m) = |x - x_U|$ and
  $n_\m = \frac{x - x_U}{|x - x_U|}$.
  
  As $\rho_*$ is the inverse of $\theta_*$, we get
  $\theta_m = \theta_*(|x - x_U|)$ on $U$, meaning
  $m = m_*(\cdot - x_U)$ on $U$.  By continuity, we even get
  $m = m_*(\cdot - x_U)$ on $\partial U \subset \{m = \pm e_3\}$.  The
  map $m_*(\cdot - x_U)$ takes the values $\pm e_3$ only at $x_U$, so
  that $\partial U \subset \{x_U\}$.  But this means that
  $U = \R^2 \setminus \{x_U\}$ is the unique connected component of
  $\Omega_m$ and thus $m = m_*(\cdot - x_0)$ on $\R^2$ with
  $x_0 := x_U$.  
\end{proof}

\section*{Declarations}

\noindent
\textit{Data Availability:}
The manuscript contains no associated data.
\par\smallskip

\noindent
\textit{Conflict of Interest:}
The authors declare that they have no conflicts of interest.
\par\smallskip

\noindent
\textit{Author Contributions:}
All authors contributed equally to this work.
\par\smallskip


\bibliography{../../nonlin,../../mura,../../stat}

\bibliographystyle{plain}

\end{document}